\documentclass{article}

\usepackage[table,dvipsnames]{xcolor}

\usepackage{lmodern}

\usepackage[margin=2.5cm]{geometry}
\usepackage[utf8]{inputenc}
\usepackage[T1]{fontenc}
\usepackage{scrbase}
\usepackage{microtype}
\usepackage[shortlabels]{enumitem}
\usepackage{etoolbox}
\usepackage{graphicx}
\usepackage{tabularx}
\usepackage{booktabs}
\usepackage{arydshln}
\usepackage{multirow}
\usepackage{makecell}
\usepackage{pdflscape}
\usepackage{authblk}

\usepackage[
backend=biber,
style=alphabetic,
maxnames=99,
backref=true
]{biblatex}
\usepackage{amsmath}
\usepackage{amssymb}
\usepackage{amsthm}
\usepackage{thmtools}
\usepackage{mathtools}
\usepackage{braket}

\usepackage{todonotes}
\usepackage{lipsum}

\usepackage{array}
\newcolumntype{Y}{>{\centering\arraybackslash}X}

\definecolor{betterblue}{RGB}{0, 110, 230}
\definecolor{bettergreen}{RGB}{164, 180, 43}

\usepackage[most]{tcolorbox}
\newtcolorbox{mybox}{
    colback=orange!10, colframe=orange!40, arc=4mm, boxrule=1pt,
    left=6pt, right=6pt, top=6pt, bottom=6pt,
    width=0.95\textwidth, center
}

\usepackage{hyperref}
\hypersetup{%
	pdftoolbar=false,
	pdfmenubar=false,
	colorlinks=true,
	urlcolor={OliveGreen},
	linkcolor={betterblue},
	citecolor={orange}
}
\usepackage[capitalise,noabbrev]{cleveref}

\theoremstyle{plain}
\newtheorem{theorem}{Theorem}[section]
\newtheorem{lemma}[theorem]{Lemma}
\newtheorem{proposition}[theorem]{Proposition}
\newtheorem{corollary}[theorem]{Corollary}

\newtheorem*{question*}{Question}
\newtheorem*{conjecture*}{Conjecture}

\AtEndEnvironment{proposition}{\setcounter{claim}{0}}
\AtEndEnvironment{lemma}{\setcounter{claim}{0}}
\AtEndEnvironment{theorem}{\setcounter{claim}{0}}

\theoremstyle{definition}
\newtheorem{definition}[theorem]{Definition}

\newtheorem{remark}[theorem]{Remark}
\newtheorem{examplex}[theorem]{Example}
\newenvironment{example}
  {\pushQED{\qed}\examplex}
  {\popQED\endexamplex}

\renewcommand{\epsilon}{\varepsilon}

\newcommand{\C}{\ensuremath{\mathbb{C}}}
\newcommand{\R}{\ensuremath{\mathbb{R}}}

\newcommand{\FdHilb}{\mathbf{FdHilb}}
\newcommand{\TFdHilb}{\mathbf{2FdHilb}}

\DeclareMathOperator{\id}{id}
\DeclareMathOperator{\Tr}{Tr}
\DeclareMathOperator{\diag}{diag}
\DeclareMathOperator{\lspan}{span}

\DeclareMathOperator{\img}{img}
\DeclareMathOperator{\rk}{rk}

\DeclareMathOperator{\MinCut}{MinCut}
\DeclareMathOperator{\MaxFlow}{MaxFlow}
\DeclareMathOperator{\GenMinCut}{MinGenCut}
\DeclareMathOperator{\val}{val}
\DeclareMathOperator{\cut}{cut}

\newcommand{\transpose}{^\mathsf{T}}
\newcommand{\adjoint}{^\dagger}
\newcommand{\inverse}{^{-1}}
\newcommand{\identity}{I}
\newcommand{\I}{\identity}
\newcommand{\compl}[1]{#1^c}
\newcommand{\op}{^{\operatorname{op}}}

\newcommand{\IFF}{\Longleftrightarrow}

\newcommand{\abs}[1]{\lvert #1 \rvert}
\newcommand{\norm}[1]{\lVert #1 \rVert}
\newcommand{\Cstar}{C\ensuremath{^*}}

\newcommand{\Hilb}{\mathcal{H}}
\newcommand{\psubseteq}{\sqsubseteq}

\newcommand{\pproj}[1]{\widetilde{#1}}
\newcommand{\mmap}[1]{\widehat{#1}}

\newcommand{\diff}{\mathop{}\!\mathrm{d}}

\newcommand{\Gsym}{G^\mathrm{sym}}
\newcommand{\Gasym}{G^\mathrm{asym}}

\newcommand{\ccdot}{\makebox[1ex]{\textbf{$\cdot$}}}

\newcommand{\lrangle}[1]{\langle #1 \rangle}

\makeatletter
\newsavebox{\bR@D}\newsavebox{\bR@T}\newsavebox{\bR@S}\newsavebox{\bR@SS}
\newcommand{\bR@tuck}{0.8}
\newcommand{\bR@raw}[1]{%
  \begingroup
  \tikzset{every picture/.style={},every node/.style={}}%
  \tikz[baseline=(n.base)]{%
    \node[inner sep=0pt,outer sep=0pt,cm={-1,0,0.25,1,(0,0)},transform shape]
      (n) {$#1\mathrm{R}$};}%
  \endgroup}
\AtBeginDocument{%
  \sbox\bR@D{\bR@raw\displaystyle}%
  \sbox\bR@T{\bR@raw\textstyle}%
  \sbox\bR@S{\bR@raw\scriptstyle}%
  \sbox\bR@SS{\bR@raw\scriptscriptstyle}}
\newcommand{\bR@use}[1]{\makebox[\bR@tuck\wd#1][l]{\usebox#1}}
\newcommand{\bR}{{\mathord{\mathchoice
{\bR@use\bR@D}{\bR@use\bR@T}{\bR@use\bR@S}{\bR@use\bR@SS}}}}
\makeatother

\makeatletter
\let\amsmath@bigm\bigm

\renewcommand{\bigm}[1]{%
  \ifcsname fenced@\string#1\endcsname
    \expandafter\@firstoftwo
  \else
    \expandafter\@secondoftwo
  \fi
  {\expandafter\amsmath@bigm\csname fenced@\string#1\endcsname}%
  {\amsmath@bigm#1}%
}

\newcommand{\DeclareFence}[2]{\@namedef{fenced@\string#1}{#2}}
\makeatother

\DeclareFence{\mid}{|} 

\makeatletter
\DeclareRobustCommand\ket[1]{%
  \@ifnextchar\bra{\k@t{#1}\!}{\k@t{#1}}%
}
\newcommand\k@t[1]{{|{#1}\rangle}}
\makeatother

\makeatletter
\newcommand{\vast}{\bBigg@{3}}
\newcommand{\Vast}{\bBigg@{4}}
\newcommand{\VVast}{\bBigg@{5}}
\newcommand{\VAST}{\bBigg@{6}}
\makeatother

\makeatletter
\patchcmd{\@maketitle}{\LARGE \@title}{\fontsize{19}{21.2}\selectfont\@title}{}{}
\makeatother

\usetikzlibrary{patterns, patterns.meta}

\input{./style/my_tikzit.sty}

\tikzstyle{myBlue}=[draw={rgb,255: red,25; green,132; blue,232}]
\tikzstyle{myBlueFill}=[fill={rgb,255: red,25; green,132; blue,232}]
\tikzstyle{myGreen}=[draw={rgb,255: red,164; green,180; blue,43}]
\tikzstyle{myPink}=[draw={rgb,255: red,208; green,12; blue,133}]
\tikzstyle{myPinkFill}=[fill={rgb,255: red,208; green,12; blue,133}]
\tikzstyle{myYellow}=[draw={rgb,255: red,255; green,160; blue,0}]
\tikzstyle{myYellowFill}=[fill={rgb,255: red,255; green,160; blue,0}]
\tikzstyle{myBrown}=[draw={rgb,255: red,173; green,100; blue,82}]
\tikzstyle{myBrownFill}=[fill={rgb,255: red,173; green,100; blue,82}]
\tikzstyle{myOrange}=[draw={rgb,255: red,245; green,124; blue,0}]
\tikzstyle{myOrangeFill}=[fill={rgb,255: red,245; green,124; blue,0}]
\tikzstyle{myPurple}=[draw={rgb,255: red,178; green,121; blue,219}]
\tikzstyle{myPurpleFill}=[fill={rgb,255: red,178; green,121; blue,219}]
\tikzstyle{myPlum}=[draw={rgb,255: red,123; green,31; blue,162}]
\tikzstyle{myPlumFill}=[fill={rgb,255: red,123; green,31; blue,162}]
\tikzstyle{myTeal}=[draw={rgb,255: red,0; green,151; blue,167}]
\tikzstyle{map}=[myBlue, fill=white, shape=circle, thick, tikzit category=maps, minimum size=20pt, inner sep=1pt, tikzit draw={rgb,255: red,25; green,132; blue,232}]
\tikzstyle{wmap}=[myBlue, fill=white, shape=rectangle, rounded corners=10pt, thick, tikzit category=maps, minimum size=20pt, inner sep=1pt, tikzit draw={rgb,255: red,25; green,132; blue,232}]
\tikzstyle{correlation}=[map, myPink, tikzit draw={rgb,255: red,208; green,12; blue,133}, tikzit category=games]
\tikzstyle{game}=[map, myPurple, tikzit draw={rgb,255: red,104; green,159; blue,56}, tikzit category=games]
\tikzstyle{q-func}=[map, myPink, tikzit draw={rgb,255: red,208; green,12; blue,133}, tikzit category=games]
\tikzstyle{frob}=[fill=white, tikzit draw=black, shape=circle, very thick, tikzit category=frobenius, minimum size=4pt, inner sep=0pt]
\tikzstyle{x-frob}=[frob, myYellow, tikzit draw={rgb,255: red,255; green,160; blue,0}, tikzit category=frobenius]
\tikzstyle{dark-x-frob}=[frob, myYellow, myYellowFill, tikzit fill={rgb,255: red,255; green,160; blue,0}, tikzit draw={rgb,255: red,255; green,160; blue,0}, tikzit category=frobenius]
\tikzstyle{y-frob}=[frob, myBrown, tikzit draw={rgb,255: red,173; green,100; blue,82}, tikzit category=frobenius]
\tikzstyle{xy-frob}=[frob, myOrange, ultra thick, tikzit draw={rgb,255: red,245; green,124; blue,0}, tikzit category=frobenius]
\tikzstyle{half-xy-proj}=[xy-frob, shape=rectangle, minimum size=6pt, tikzit draw={rgb,255: red,245; green,124; blue,0}, tikzit category=frobenius, tikzit fill=white, tikzit shape=rectangle]
\tikzstyle{comp-half-xy-proj}=[xy-frob, shape=rectangle, minimum size=6pt, tikzit draw={rgb,255: red,245; green,124; blue,0}, tikzit category=frobenius, diagonal strike, tikzit fill=white, tikzit shape=rectangle]
\tikzstyle{half-xy-proj-flipped}=[half-xy-proj, myOrangeFill, tikzit draw={rgb,255: red,245; green,124; blue,0}, tikzit category=frobenius, tikzit fill={rgb,255: red,245; green,124; blue,0}, tikzit shape=rectangle]
\tikzstyle{twisted-xy-proj}=[myOrange, left color=white, right color={rgb,255: red,245; green,124; blue,0}, tikzit draw={rgb,255: red,245; green,124; blue,0}, shape=rectangle, thin, tikzit category=frobenius, minimum size=8pt, inner sep=0pt, tikzit fill=white, tikzit shape=rectangle]
\tikzstyle{twisted-xy-proj-flipped}=[myOrange, right color=white, left color={rgb,255: red,245; green,124; blue,0}, tikzit draw={rgb,255: red,245; green,124; blue,0}, shape=rectangle, thin, tikzit category=frobenius, minimum size=8pt, inner sep=0pt, tikzit fill={rgb,255: red,245; green,124; blue,0}, tikzit shape=rectangle]
\tikzstyle{a-frob}=[frob, myPurple, tikzit draw={rgb,255: red,178; green,121; blue,219}, tikzit category=frobenius]
\tikzstyle{dark-a-frob}=[frob, myPurple, myPurpleFill, tikzit fill={rgb,255: red,178; green,121; blue,219}, tikzit draw={rgb,255: red,178; green,121; blue,219}, tikzit category=frobenius]
\tikzstyle{b-frob}=[frob, myBlue, tikzit draw={rgb,255: red,25; green,118; blue,210}, tikzit category=frobenius]
\tikzstyle{ab-frob}=[frob, myPlum, ultra thick, tikzit draw={rgb,255: red,123; green,31; blue,162}, tikzit category=frobenius]
\tikzstyle{half-ab-proj}=[ab-frob, shape=rectangle, minimum size=6pt, tikzit draw={rgb,255: red,123; green,31; blue,162}, tikzit category=frobenius, tikzit shape=rectangle]
\tikzstyle{comp-half-ab-proj}=[ab-frob, shape=rectangle, diagonal strike, minimum size=6pt, tikzit draw={rgb,255: red,123; green,31; blue,162}, tikzit category=frobenius, tikzit shape=rectangle]
\tikzstyle{half-ab-proj-flipped}=[half-ab-proj, myPlumFill, tikzit category=frobenius, tikzit draw={rgb,255: red,123; green,31; blue,162}, tikzit fill={rgb,255: red,123; green,31; blue,162}, tikzit shape=rectangle]
\tikzstyle{twisted-ab-proj}=[myPlum, left color=white, right color={rgb,255: red,123; green,31; blue,162}, tikzit draw={rgb,255: red,123; green,31; blue,162}, shape=rectangle, thin, tikzit category=frobenius, minimum size=8pt, inner sep=0pt, tikzit fill=white, tikzit shape=rectangle]
\tikzstyle{twisted-ab-proj-flipped}=[myPlum, right color=white, left color={rgb,255: red,123; green,31; blue,162}, shape=rectangle, thin, tikzit category=frobenius, minimum size=8pt, inner sep=0pt, tikzit draw={rgb,255: red,123; green,31; blue,162}, tikzit fill={rgb,255: red,123; green,31; blue,162}, tikzit shape=rectangle]
\tikzstyle{trap}=[fill=white, myPink, shape=trapezium, thick, tikzit shape=rectangle, trapezium left angle=90, trapezium right angle=70, minimum height=16pt, minimum width=16pt, inner sep=2pt, draw={rgb,255: red,25; green,132; blue,232}, tikzit category=maps]
\tikzstyle{conj}=[trap, trapezium left angle=70, trapezium right angle=90, tikzit draw={rgb,255: red,25; green,132; blue,232}, tikzit shape=rectangle, tikzit category=maps]
\tikzstyle{dagger}=[conj, shape border rotate=180, tikzit draw={rgb,255: red,25; green,132; blue,232}, tikzit shape=rectangle, tikzit category=maps]
\tikzstyle{transpose}=[trap, shape border rotate=180, tikzit draw={rgb,255: red,25; green,132; blue,232}, tikzit shape=rectangle, tikzit category=maps]
\tikzstyle{strap}=[fill=white, myPink, shape=trapezium, thick, tikzit shape=rectangle, trapezium left angle=90, trapezium right angle=70, draw={rgb,255: red,25; green,132; blue,232}, tikzit category=maps]
\tikzstyle{sconj}=[strap, trapezium left angle=70, trapezium right angle=90, tikzit draw={rgb,255: red,25; green,132; blue,232}, tikzit shape=rectangle, tikzit category=maps]
\tikzstyle{sdagger}=[sconj, shape border rotate=180, tikzit draw={rgb,255: red,25; green,132; blue,232}, tikzit shape=rectangle, tikzit category=maps]
\tikzstyle{stranspose}=[strap, shape border rotate=180, tikzit draw={rgb,255: red,25; green,132; blue,232}, tikzit shape=rectangle, tikzit category=maps]
\tikzstyle{wbra}=[myTeal, fill=white, tikzit draw={rgb,255: red,0; green,151; blue,167}, isosceles triangle,
    isosceles triangle apex angle=95, tikzit category=maps, shape border rotate=90, inner sep=0pt, minimum size=20pt]
\tikzstyle{bra}=[myTeal, fill=white, tikzit draw={rgb,255: red,0; green,151; blue,167}, regular polygon, regular polygon sides=3, tikzit category=maps, inner sep=0pt, minimum size=20pt]
\tikzstyle{bra_s}=[myTeal, fill=white, tikzit draw={rgb,255: red,0; green,151; blue,167}, regular polygon, regular polygon sides=3, tikzit category=maps, inner sep=-1.25pt, minimum size=20pt]
\tikzstyle{bra_xs}=[myTeal, fill=white, tikzit draw={rgb,255: red,0; green,151; blue,167}, regular polygon, regular polygon sides=3, tikzit category=maps, inner sep=-1.75pt, minimum size=20pt]
\tikzstyle{bra_xxs}=[myTeal, fill=white, tikzit draw={rgb,255: red,0; green,151; blue,167}, regular polygon, regular polygon sides=3, tikzit category=maps, inner sep=-3pt, minimum size=1pt]
\tikzstyle{ket}=[bra, shape border rotate=180, tikzit draw={rgb,255: red,0; green,151; blue,167}, tikzit category=maps]
\tikzstyle{wket}=[wbra, shape border rotate=270, tikzit draw={rgb,255: red,0; green,151; blue,167}, tikzit category=maps]
\tikzstyle{ket_s}=[{bra_s}, shape border rotate=180, tikzit draw={rgb,255: red,0; green,151; blue,167}, tikzit category=maps]
\tikzstyle{ket_xs}=[{bra_xs}, shape border rotate=180, tikzit draw={rgb,255: red,0; green,151; blue,167}, tikzit category=maps]
\tikzstyle{ket_xxs}=[{bra_xxs}, shape border rotate=180, tikzit draw={rgb,255: red,0; green,151; blue,167}, tikzit category=maps]
\tikzstyle{map_multi}=[myGreen, thick, shape=rectangle, tikzit draw={rgb,255: red,25; green,132; blue,232}, tikzit shape=rectangle, minimum width=1.5cm, tikzit category=maps]
\tikzstyle{small_map_multi}=[{map_multi}, minimum width=1cm, tikzit category=maps, shape=rectangle, tikzit draw={rgb,255: red,25; green,132; blue,232}]
\tikzstyle{proj}=[fill=white, shape=rectangle, tikzit category=maps, myGreen]
\tikzstyle{cket}=[bra, shape=cut triangle down, tikzit draw={rgb,255: red,0; green,128; blue,128}, apex angle=110, minimum width=.9cm]
\tikzstyle{cbra}=[cket, tikzit draw={rgb,255: red,0; green,128; blue,128}, shape=cut triangle up]
\tikzstyle{wcket}=[cket, apex angle=135, vertical offset=.8, minimum width=1.6cm]
\tikzstyle{wcbra}=[cbra, apex angle=135, vertical offset=0, minimum width=1.6cm]
\tikzstyle{ccket}=[cket, shape=cut triangle down alt, tikzit draw={rgb,255: red,0; green,128; blue,128}, tikzit shape=circle]
\tikzstyle{ccbra}=[cbra, shape=cut triangle up alt, tikzit draw={rgb,255: red,0; green,128; blue,128}, tikzit shape=circle]
\tikzstyle{idendo}=[fill=white, draw={black!80}, shape=circle, dashed, inner sep=0.1]
\tikzstyle{selfadj}=[trap, signal, signal pointer angle=120, signal to=east, fill=white, tikzit shape=rectangle, draw={rgb,255: red,25; green,132; blue,232}]
\tikzstyle{selfadj-conj}=[selfadj, signal to=west, fill=white, tikzit draw=blue, tikzit shape=rectangle]
\tikzstyle{pselfadj}=[selfadj, myGreen, tikzit draw={rgb,255: red,64; green,161; blue,0}, tikzit shape=rectangle]
\tikzstyle{pconjselfadj}=[selfadj, signal to=west, myGreen, tikzit draw={rgb,255: red,64; green,161; blue,0}, tikzit shape=rectangle]
\tikzstyle{pselfadjselfdual}=[tikzit draw={rgb,255: red,64; green,161; blue,0}, signal, signal pointer angle=120, signal to=east and west, fill=white, tikzit shape=circle, myGreen]
\tikzstyle{x-action}=[myOrange, fill=white, shape=rectangle, tikzit draw={rgb,255: red,245; green,124; blue,0}]
\tikzstyle{a-action}=[myPurple, fill=white, shape=rectangle, tikzit draw={rgb,255: red,178; green,121; blue,219}]
\tikzstyle{y-action}=[myBrown, fill=white, shape=rectangle, tikzit draw={rgb,255: red,173; green,100; blue,82}]

\tikzstyle{makesgaps}=[preaction={draw=white, line width=2.5pt}]
\tikzstyle{makeslargegaps}=[preaction={draw=white, line width=3.5pt}]
\tikzstyle{hilbert}=[-, myTeal, thick, tikzit draw={rgb,255: red,0; green,151; blue,167}]
\tikzstyle{hilbert2}=[-, myGreen, thick, tikzit draw={rgb,255: red,148; green,210; blue,75}]
\tikzstyle{x-wire}=[-, myYellow, very thick, tikzit draw={rgb,255: red,255; green,160; blue,0}]
\tikzstyle{x-gap-wire}=[-, makesgaps, myYellow, very thick, tikzit draw={rgb,255: red,255; green,160; blue,0}]
\tikzstyle{y-wire}=[-, myBrown, very thick, tikzit draw={rgb,255: red,173; green,100; blue,82}]
\tikzstyle{xy-wire}=[-, myOrange, line width=2pt, tikzit draw={rgb,255: red,245; green,124; blue,0}]
\tikzstyle{xy-gap-wire}=[-, makeslargegaps, myOrange, line width=2pt, tikzit draw={rgb,255: red,245; green,124; blue,0}]
\tikzstyle{a-wire}=[-, myPurple, very thick, tikzit draw={rgb,255: red,178; green,121; blue,219}]
\tikzstyle{a-gap-wire}=[-, makesgaps, myPurple, very thick, tikzit draw={rgb,255: red,178; green,121; blue,219}]
\tikzstyle{b-wire}=[-, myBlue, very thick, tikzit draw={rgb,255: red,25; green,118; blue,210}]
\tikzstyle{ab-wire}=[-, myPlum, tikzit draw={rgb,255: red,123; green,31; blue,162}, line width=2pt]
\tikzstyle{ab-gap-wire}=[-, makeslargegaps, myPlum, tikzit draw={rgb,255: red,123; green,31; blue,162}, line width=2pt]
\tikzstyle{hilbertnoarrow}=[-, tikzit draw={rgb,255: red,0; green,151; blue,167}, thick, myTeal]
\tikzstyle{dirwire}=[-]
\tikzstyle{wire}=[-]
\tikzstyle{dirwire02}=[-, tikzit draw=magenta]
\tikzstyle{invis}=[-, draw=none, tikzit draw={rgb,255: red,128; green,128; blue,128}]
\tikzstyle{z-wire}=[-, dashed, draw={rgb,255: red,64; green,64; blue,64}, very thick]
\tikzstyle{dir-a-wire}=[-, a-wire, tikzit draw={rgb,255: red,178; green,121; blue,219}]
\tikzstyle{dir-b-wire}=[-, b-wire, tikzit draw={rgb,255: red,25; green,118; blue,210}]
\tikzstyle{dir-x-wire}=[-, x-wire, tikzit draw={rgb,255: red,255; green,160; blue,0}]
\tikzstyle{dir-y-wire}=[-, y-wire, tikzit draw={rgb,255: red,173; green,100; blue,82}]
\tikzstyle{f-dir-a-wire}=[-, a-wire, tikzit draw={rgb,255: red,178; green,121; blue,219}]
\tikzstyle{f-dir-x-wire}=[-, x-wire, tikzit draw={rgb,255: red,255; green,160; blue,0}]
\tikzstyle{f-dir-y-wire}=[-, y-wire, tikzit draw={rgb,255: red,173; green,100; blue,82}]
\tikzstyle{b-dir-a-wire}=[-, a-wire, tikzit draw={rgb,255: red,178; green,121; blue,219}]
\tikzstyle{b-dir-x-wire}=[-, x-wire, tikzit draw={rgb,255: red,255; green,160; blue,0}]
\tikzstyle{b-dir-y-wire}=[-, y-wire, tikzit draw={rgb,255: red,173; green,100; blue,82}]

\input{./style/tikz_adjustments.sty}

\input{./style/shapes.tex}

\renewcommand{\emptyset}{\varnothing}

\usepackage{abstract}

\title{A Max-Flow Min-Cut Theorem for Quantum Graphs}

\author[1,2]{Gian Luca Spitzer \thanks{\href{mailto:gian-luca.spitzer@u-bordeaux.fr}{gian-luca.spitzer@u-bordeaux.fr}}}
\author[2]{Ion Nechita}
\affil[1]{LaBRI, Universit\'e de Bordeaux, CNRS, Bordeaux INP, UMR-5800, France}
\affil[2]{Laboratoire de Physique Th\'eorique, Universit\'e de Toulouse, CNRS, UPS, France}

\date{}

\begin{document}

\maketitle

\begin{abstract}
  We propose natural definitions of cuts and flows for quantum graphs. For cuts
  we consider two notions: ordinary cuts, which partition the vertex set into
  two parts; and \emph{generalised} cuts, defined as a continuous relaxation of
  the binary partition function. We show that our definitions satisfy a
  max-flow min-cut theorem: For any disjoint terminal sets $S$ and $T$, the maximum value of an $S$-$T$-flow of a quantum
  graph equals the minimum value of a generalised $S$-$T$ cut, while the
  minimum value of an ordinary $S$-$T$ cut can be larger. While for classical
  graphs the values of the minimum generalised cut and the minimum ordinary cut
  coincide---recovering the classical max-flow min-cut theorem---we also show how it
  is more natural to think of flows as duals of generalised cuts instead of
  ordinary cuts. 
\end{abstract}

\tableofcontents

\section{Introduction}\label{sec:intro}
Since their introduction in the context of quantum Shannon theory
\cite{duan_zeroerror_2013} and independently in the context of quantum metric
spaces \cite{weaver_quantum_2010}, quantum graphs have become an active area of
research. Naturally, given their classical counterparts, much of this research
has focused on generalising concepts of classical graph theory to the quantum
case. Examples include basic notions such as connectivity
\cite{chavezdominguez_connectivity_2021,courtney_connectivity_2026} and
colourings \cite{brannan_quantumclassical_2022,ganesan_spectral_2023}; more
complex ideas such as Cayley graphs \cite{wasilewski_quantum_2024} and graph
non-local games \cite{brannan_quantumclassical_2022,goldberg_quantum_2026}; and
even interdisciplinary topics such as graph algebras
\cite{brannan_quantum_2023} and association schemes
\cite{gromada_quantum_2024}. Famous theorems from classical
graph theory have also been generalised, for instance Frucht's theorem
\cite{brannan_quantum_2025} or Ramsey's theorem \cite{weaver_quantum_2017}.

In this work, we extend this landscape by studying \emph{cuts} and \emph{flows} of quantum
graphs. Classically, a \emph{flow network} is a graph whose edges are assigned
capacities. Given a set of source vertices and a set of sink vertices, the
goal is to assign to each edge a flow in order to maximise the total flow from
the source set to the sink set. The assigned flows must respect the edge capacities,
and at each vertex that is not a source or a sink the total incoming flow must
equal the outgoing flow. Clearly, flow networks can be used to model a large
array of real-world problems such as network routing, logistics optimisation,
and scheduling problems \cite{ahuja1993network,ford1962flows}.
Cuts, on the other hand, are about a failure mode of these problems: They
identify bottlenecks. A cut is a set of edges that disconnects source from sink
when removed. The cost of a cut is then given by the total capacity of the removed
edges. Clearly, any flow is upper bounded by the cost of any cut and hence by
the cost of the smallest cut. The famous \emph{Max-Flow Min-Cut theorem}
\cite{ford1956maximal,dantzig2003max} tells us that this upper bound can always be achieved.

Our approach to generalising the classical framework of cuts and flows to the
quantum setting centers on the following four concepts.
\begin{enumerate}
  \item \textbf{Terminals.} A fundamental requirement for both cuts and flows
    is a notion of source and sink sets. In the quantum setting, their role is
    played by projective subsets of the quantum vertex set, as introduced in
    \cite{spitzer_unification_2026}.
  \item \textbf{Capacity Maps.} Completely positive symmetric operators $K
    \colon X \to X$. These are the quantum analogue of capacity matrices. The
    complete positivity captures that all capacities are non-negative while the
    symmetry captures that capacities are assigned to undirected edges. 
  \item \textbf{Flow Maps.} Real, antisymmetric operators $F \colon X \to X$
    satisfying a conservation condition, such that $K - F$ is completely
    positive. The antisymmetry captures that flows are directed, and that
    traversing an edge in the opposite direction negates the assigned flow. The
    complete positivity constraint encodes that flows cannot exceed the edge
    capacities.
  \item \textbf{Vertex Partitions.} Binary partitions of the quantum vertex set in
    terms of projective subsets $X = E \sqcup F$ with $E, F \sqsubseteq X$. Any such partition induces a cut by
    considering the edges between the two blocks of the partition.
\end{enumerate}
We also define a weaker notion of cuts, called \emph{generalised cuts}.
Classically, these correspond to a continuous relaxation of the partition
function defining an ordinary cut. To that end, we introduce a quantum analogue
of real-valued functions on a set. With these definitions in hand, we can then
naturally define the quantities $\MaxFlow(G, K; S, T)$, denoting the maximum
achievable flow from $S$ to $T$ on a quantum graph $G$ with capacities given by
$K$; $\MinCut(G, K; S, T)$, denoting the corresponding minimum achievable cut;
and $\GenMinCut(G, K; S, T)$, its continuous relaxation.

What emerges is an elegant theory naturally extending the classical landscape.
While for classical graphs the abovementioned quantities are naturally phrased
in terms of linear programs, in the quantum case they are described by
semidefinite programs. This perspective then allows us to show strong duality
results between flows and generalised cuts, culminating in the following
Max-Flow Min-Cut theorem for quantum graphs.

\begin{restatable*}[Quantum Max-Flow Min-Cut]{theorem}{thmmfmc}\label{thm:qmaxflowmincut}
  Let $G$ be a quantum graph. Then for all capacities $K$ on $G$ and non-trivial disjoint $S, T \psubseteq V(G)$ it holds that
  \begin{equation*}
    \MaxFlow(G, K; S, T) = \GenMinCut(G, K; S, T) \leq \MinCut(G, K; S, T).
  \end{equation*}
\end{restatable*}

The distinction between the minimum cut and the minimum generalised cut is
particularly interesting. It follows from our theorem together with the classical
max-flow min-cut theorem that these quantities are equal for classical graphs.
We show that this is generally not the case quantumly. This is because capacity
maps may enforce coherences between different edges that ordinary cuts are
blind to, while generalised cuts may exploit them.

\smallskip

We conclude by computing flows and cuts for a number of natural examples. These
include the complete quantum graphs on arbitrary quantum vertex sets, and a
subset of the quantum graphs having classical symmetry studied in \cite{spitzer2026quantum}. We show that
for certain classes of quantum graphs equality between minimum cuts and minimum
generalised cuts can be recovered. On the other hand, we show that our
framework exhibits non-classical phenomena, such as the existence of connected
quantum graphs admitting no non-zero flows for certain choices of source and
sink sets.

\medskip

The paper is organised as follows. In \cref{sec:prelims}, we recall the most
important ideas underlying our treatment of quantum graphs. This includes a
detailed treatment of the realignment procedure that was already used in
\cite{spitzer2026quantum} to classify quantum graphs invariant under classical
group actions. We show why the same concept can take different forms and how
they arise. In \cref{sec:cutsflows}, we introduce our definitions of capacity
maps, flows, and cuts and derive fundamental results. Centrally, in
\cref{sec:mncase} we show how the quantum theory of cuts and flows mirrors the
classical landscape of representing these problems as linear programs,
culminating in the statement and proof of the max-flow min-cut theorem for
quantum graphs. In \cref{sec:examples}, we compute flows and cuts of example
quantum graphs. Finally, \cref{sec:conclusion} summarises our findings and
discusses potential avenues for future research.

\section{Preliminaries}\label{sec:prelims}
First, we fix some general conventions. We assume all projectors to be orthogonal projectors: A linear map $P \colon H \to H$ on a Hilbert space $H$ is a \emph{projector} if $P^2 = P\adjoint = P$. Similarly, we call an element $p$ of a \Cstar-algebra a \emph{projection} if $p^2 = p^* = p$. When talking about classical graphs, we assume the vertex set to be $[n] \coloneqq \{1, \dots, n\}$. We also assume all classical graphs to be undirected and loopless. We write $i \sim j$ if the vertices $i$ and $j$ are adjacent in the graph.

\subsection{The Graphical Formalism}

We will work in the diagrammatic framework of quantum sets and quantum functions introduced by Musto, Reutter, and Verdon \cite{musto_compositional_2018}. We recall the main concepts here briefly. For a more complete treatment we refer to their paper, and for the full category-theoretic background to \cite{vicary_categorical_2011}. The main idea is to treat quantum sets as special objects in the category $\FdHilb$.

\begin{definition}\label{def:qset}
  A \emph{quantum set} is a special symmetric $\dagger$-Frobenius monoid in $\FdHilb$.
\end{definition}

As a $\dagger$-compact category, $\FdHilb$ comes equipped with a string-diagrammatic calculus that can be used in the place of symbolic equations to derive statements about the category. The primitives of the calculus are \emph{strings}, which represent objects, and \emph{boxes}, which represent morphisms. 
\begin{equation*}
  \begin{tikzpicture}
	\begin{pgfonlayer}{nodelayer}
		\node [style=none] (0) at (0, -0.75) {};
		\node [style=none] (2) at (0, 0.75) {};
	\end{pgfonlayer}
	\begin{pgfonlayer}{edgelayer}
		\draw [style=x-wire] (0.center) to (2.center);
	\end{pgfonlayer}
\end{tikzpicture}
 \qquad\qquad \begin{tikzpicture}
	\begin{pgfonlayer}{nodelayer}
		\node [style=map] (0) at (0, 0) {$f$};
		\node [style=none] (1) at (0, 0.75) {};
		\node [style=none] (2) at (0, -0.75) {};
	\end{pgfonlayer}
	\begin{pgfonlayer}{edgelayer}
		\draw [style=x-wire] (2.center) to (0);
		\draw [style=a-wire] (0) to (1.center);
	\end{pgfonlayer}
\end{tikzpicture}

\end{equation*}
We use special notation for the multiplication and unit morphisms of monoids, which we represent by fusion of two strings and dots respectively.
\begin{equation*}
  \input{diagrams/mult.tikz} \qquad\qquad \begin{tikzpicture}
	\begin{pgfonlayer}{nodelayer}
		\node [style=x-frob] (0) at (0, -0.25) {};
		\node [style=none] (2) at (0, 0.25) {};
	\end{pgfonlayer}
	\begin{pgfonlayer}{edgelayer}
		\draw [style=x-wire] (0) to (2.center);
	\end{pgfonlayer}
\end{tikzpicture}

\end{equation*}
Composition $g \circ f$ of functions is represented by serial composition of diagrams, while tensor products $f \otimes h$ are represented by parallel composition. Diagrams are read from bottom to top.
\begin{equation*}
  \input{diagrams/composition.tikz} \qquad\qquad \input{diagrams/tensorprod.tikz}
\end{equation*}
The monoidal unit, $\C$ in the case of $\FdHilb$, is represented by the empty diagram.

Crucially, any equation that can be derived in the string-diagrammatic calculus by moving and bending strings, moving boxes and sliding them along strings, and rotating boxes also holds in $\FdHilb$, see \cite[Theorem 3.28]{heunen_categories_2019}. The action of rotating boxes by $180^\circ$ corresponds to taking the dual of the corresponding morphism. Since we work in a category of Hilbert spaces and linear maps, we refer to this operation as the \emph{transpose} and denote it by $(-)\transpose$. Similarly, we may mirror a diagram across a horizontal axis, which corresponds to application of the $\dagger$-functor to the corresponding morphism. On the level of linear maps, this is the Hilbert space \emph{adjoint}. Finally, we can compose the two functors. This yields the \emph{conjugate functor} $\overline{(-)} = ((-)\adjoint)\transpose = ((-)\transpose)\adjoint$. Graphically, this action must then correspond to mirroring across a vertical axis. In the spirit of these geometric representations, we will sometimes denote morphisms by asymmetric boxes and write
\begin{equation*}
  \input{diagrams/nonsym-func.tikz}~, \qquad \input{diagrams/nonsym-func-transpose.tikz}~, \qquad \input{diagrams/nonsym-func-adjoint.tikz}~, \qquad \input{diagrams/nonsym-func-conj.tikz}~.
\end{equation*}
Invariance of a morphism under any of these operations can then be encoded in the symmetries of the box representing it. For example, a self-adjoint morphism may be represented as follows.
\begin{equation*}
  \begin{tikzpicture}
	\begin{pgfonlayer}{nodelayer}
		\node [style=selfadj] (0) at (0, 0) {$f$};
		\node [style=none] (1) at (0, 0.75) {};
		\node [style=none] (2) at (0, -0.75) {};
	\end{pgfonlayer}
	\begin{pgfonlayer}{edgelayer}
		\draw [style=x-wire] (2.center) to (0);
		\draw [style=x-wire, in=270, out=90] (0) to (1.center);
	\end{pgfonlayer}
\end{tikzpicture}

\end{equation*}
Elements of a Hilbert space $X$ can be represented in the graphical calculus by morphisms $\C \to X$. We represent them by triangles with a corner cut off, unless they are self-conjugate, in which case we represent them by a full triangle.
\begin{equation*}
  \begin{tikzpicture}
	\begin{pgfonlayer}{nodelayer}
		\node [style=cket] (0) at (0, -0.25) {$x$};
		\node [style=none] (1) at (0, 0.5) {};
	\end{pgfonlayer}
	\begin{pgfonlayer}{edgelayer}
		\draw [style=x-wire] (0) to (1.center);
	\end{pgfonlayer}
\end{tikzpicture}
 \qquad\qquad \begin{tikzpicture}
	\begin{pgfonlayer}{nodelayer}
		\node [style={ket_s}] (0) at (0, -0.25) {$\;x\;$};
		\node [style=none] (2) at (0, 0.5) {};
	\end{pgfonlayer}
	\begin{pgfonlayer}{edgelayer}
		\draw [style=x-wire] (0) to (2.center);
	\end{pgfonlayer}
\end{tikzpicture}

\end{equation*}

\begin{remark}
  A special case of diagrams are those without dangling wires. An example is the composition of an element map with its adjoint. In $\FdHilb$, these diagrams correspond to linear maps $\C \to \C$. Any such map is uniquely defined by its value at $1$, so we will view these diagrams simply as scalars.
\end{remark}

In certain cases, the graphical calculus can be refined by replacing a single wire with two parallel wires. This allows us to look ``inside of'' certain linear maps. In particular, this is the case when the object represented by the wire is $M_n$. The idea is to make use of the isomorphism $M_n \cong \C^n \otimes \C^n$ given graphically by
\begin{equation}\label{eq:mn-cncn-iso}
  M_n \to \C^n \otimes \C^n,\qquad \input{diagrams/mn-cncn-iso.tikz}\enspace.
\end{equation}
The $\dagger$-SSFM structure of $M_n$ is then inherited by $\C^n \otimes \C^n$ by defining the multiplication and unit maps as
\begin{equation*}
  \input{diagrams/mn-cncn-mult.tikz} \qquad\quad\qquad \input{diagrams/mn-cncn-unit.tikz}\enspace.
\end{equation*}
The normalisation of the multiplication map is necessary to guarantee that the resulting monoid is special, that is that $m \circ m\adjoint = \id$. The normalisation of the unit follows from the relation $m \circ (\id \otimes u) = \id$.

\begin{remark}
  The double-wire formalism highlights an important discrepancy between the categorical transpose $(-)\transpose$ and the linear algebraic transpose $(-)^t$. For morphisms $f\colon \C^n \to \C^m$ it can be verified that the categorical transpose is equal to the linear-algebraic transpose of $f$ in the standard basis. However, when computing the transpose of a morphism $f\colon M_n \to M_m$, the definitions of the duality morphisms in terms of multiplication and counit---cf. \cite[Equation 14]{musto_compositional_2018}---yield
  \begin{equation*}
    \input{diagrams/transpose.tikz}\enspace. 
  \end{equation*}
  As can be verified by plugging in the standard basis, the linear-algebraic transpose of $f$ viewed as an $m^2 \times n^2$ matrix is instead given by
  \begin{equation*}
    \input{diagrams/transpose-linalg.tikz}\enspace. 
  \end{equation*}
  Concretely, we have $f^t = \Sigma f\transpose \Sigma$, where $\Sigma$ is the flip map.
\end{remark}

\subsection{Realignment}

A particularly useful application of the double-wire calculus is the treatment of \emph{Choi matrices} \cite{jamiolkowski1972linear,choi1975completely,watrous2018theory}. Recall that given a linear operator $\Phi\colon M_n \to M_n$, its Choi matrix is defined as 
\begin{equation}\label{eq:choidef}
  J(\Phi) = \sum_{i, j = 1}^n \ket{i}\bra{j} \otimes \Phi(\ket{i}\bra{j})
\end{equation}
It turns out that we can represent the Choi matrix of $\Phi$ without sums by permuting tensor legs. Suppose that we are given a map $\Phi\colon M_n \to M_n$. Using the double-wire calculus, we can write $J(\Phi)$ graphically as
\begin{equation}\label{eq:choiraw}
  J(\Phi) \quad=\quad \sum_{i,j = 1}^n \enspace \input{diagrams/jphi-exp.tikz} \quad=\quad \sum_{i,j = 1}^n \enspace \input{diagrams/jphi-exp-2.tikz}\enspace.
\end{equation}
Now note that the unit of $\C^n$ is given by $\sum_{i = 1}^n \ket{i}$ while the comultiplication is the copy map $m\adjoint \colon \ket{i} \mapsto \ket{i} \otimes \ket{i}$. We can thus simplify \cref{eq:choiraw} as
\begin{equation*}
  \input{diagrams/jphi-exp-simpl.tikz}
\end{equation*}
where we used that comultiplying the unit corresponds to one of the duality morphisms of a $\dagger$-SSFM, see \cite[Equation 14]{musto_compositional_2018}. This gives $J(\Phi)$ as an element of $M_n \otimes M_n \cong (\C^n \otimes \C^n) \otimes (\C^n \otimes \C^n)$. It will usually be more convenient to view it as a map $\C^n \otimes \C^n \to \C^n \otimes \C^n$. We obtain such a map by unvectorising, using the isomorphism
\begin{equation*}
  \C^n \otimes \C^n \to B(\C^n),\quad \input{diagrams/unvec.tikz}\enspace.
\end{equation*}
Note that this is the inverse of (\ref{eq:mn-cncn-iso}), when viewed as a map $B(\C^n) \to \C^n \otimes \C^n$. This yields the map
\begin{equation*}
  \input{diagrams/jphi-unvec.tikz}\enspace.
\end{equation*}
For reasons that will become clear soon, we will also include a normalisation, which leads to the following definition.

\begin{definition}
  Let $\Phi\colon M_n \to M_n$ be a morphism in $\FdHilb$. Then its \emph{realignment} is given by 
  \begin{equation*}
    \input{diagrams/def-realignment.tikz}\enspace,
  \end{equation*}
  which corresponds to its normalised Choi matrix, viewed as a morphism $M_n \to M_n$. The inverse of the realignment operation is given by 
  \begin{equation*}
    \input{diagrams/def-realignment-inv.tikz}\enspace.
  \end{equation*}
\end{definition}

Note that we made two choices in our derivation of the realignment map. The first is the Choi convention of \cref{eq:choidef}. It is equally valid to apply $\Phi$ to the first tensor factor. The second convention is the choice of vectorisation map. It is equally valid to define the vectorisation map as
\begin{equation*}
  M_n \to \C^n \otimes \C^n,\qquad \input{diagrams/mn-cncn-iso-flip.tikz}\enspace.
\end{equation*}
The two maps differ by a transposition. In total, this yields four different conventions. Performing the same derivations for each leads to four different realignments. We summarise them in \cref{fig:realignmentconv}.
\begin{table}[h]
\begin{center}
\newcolumntype{C}{>{\centering\arraybackslash}m{2cm}}
  \begin{tabular}{C|C|C@{}m{0pt}@{}}
& $\id ~\otimes~ \Phi$ & $\Phi ~\otimes~ \id$ & \rule{0pt}{1cm} \\
\hline
    \input{diagrams/vec-1.tikz} & \input{diagrams/realignment-alt-4.tikz} & \input{diagrams/realignment-alt-2.tikz} & \rule{0pt}{2cm} \\
    \hline
    \input{diagrams/vec-2.tikz} & \input{diagrams/realignment-alt-3.tikz} & \input{diagrams/realignment-alt-1.tikz} & \rule{0pt}{2cm}
  \end{tabular}
\end{center} 
\caption{The different conventions for realignment procedures.}
\label{fig:realignmentconv}
\end{table}

Since the realignment is invertible, it can be seen as providing an alternative equivalent perspective on the same map. The same fundamental property may then take different forms in the two pictures. An obvious example is Choi's theorem \cite[Theorem 2]{choi1975completely}, which says that $\Phi$ is completely positive if and only if $\Phi^R$ is positive semidefinite. This correspondence can be extended to a correspondence between self-conjugate and self-adjoint maps.

\begin{proposition}\label{prop:rlgnrealselfadj}
  Let $X$ and $Y$ be $\dagger$-SSFMs and let $f\colon X \to Y$ be a morphism. Then $f$ satisfies $\overline{f} = f$ if and only if $(f^R)\adjoint = f^R$.
\end{proposition}

Another important property of the realignment is that it acts as a kind of ``Fourier transform'' between composition and the categorical Schur product. 

\begin{definition}
  Let $X$ and $Y$ be $\dagger$-SSFMs and $f, g \colon X \to Y$ be morphisms. The \emph{Schur product} $f \bullet g: X \to Y$ is defined graphically as
  \begin{equation*}
    \input{diagrams/schur-prod.tikz}\enspace.
  \end{equation*}
\end{definition}

Note that if the quantum sets are classical, this corresponds to the entrywise product of $f$ and $g$ in the standard basis.

\begin{proposition}\label{prop:rlgnschurtocomp}
  Let $X$ and $Y$ be $\dagger$-SSFMs. For all $f, g \colon X \to Y$ it holds that $(f \bullet g)^R = f^R \circ g^R$.
\end{proposition}

\cref{prop:rlgnschurtocomp} is the reason why we need to include the normalisation in the realignment procedure. Without it, the normalisation of the multiplication in $\C^n \otimes \C^n$ would lead to an additional factor $n$ on the composition side. The proofs of these results are straightforward diagrammatic derivations. Moreover, the results are independent of the realignment convention. For example, the case for the bottom left quadrant of \cref{fig:realignmentconv} was proved in \cite[Proposition 2.9]{spitzer2026quantum}.
We conclude by noting that it is also possible to define an analogous realignment procedure for maps between arbitrary $\dagger$-SSFMs beyond $M_n$. However, this requires passing to the $2$-category $\TFdHilb$ and its shaded graphical calculus (see Appendix~\ref{app:cutflowgeneral}). In this setting, Choi's theorem as well as \cref{prop:rlgnrealselfadj,prop:rlgnschurtocomp} continue to hold, see \cite{verdon_covariant_2022}. 

\subsection{Quantum Sets and Quantum Graphs}

Having introduced the necessary mathematical tools, let us now recall some results about quantum sets and quantum graphs. In \cref{def:qset}, we have defined quantum sets as $\dagger$-SSFMs in $\FdHilb$. It will sometimes be convenient to revert to the treatment of quantum sets as finite-dimensional \Cstar-algebras. The translation dictionary is quite simple and is informed by the following result by Vicary.

\begin{proposition}[Theorem 4.7 in \cite{vicary_categorical_2011}]
  The category of special symmetric $\dagger$-Frobenius monoids in $\FdHilb$ is equal to the category of finite-dimensional \Cstar-algebras. In particular, every $\dagger$-SSFM in $\FdHilb$ admits a norm that turns it into a \Cstar-algebra whose involution is represented by the conjugate functor. 
\end{proposition}

This correspondence also informs us how to represent positive elements in the graphical calculus. A positive element in a \Cstar-algebra is an element of the form $x = \sum_{i = 1}^k y_i^* y_i$. This can be represented graphically as 
\begin{equation*}
  \begin{tikzpicture}
	\begin{pgfonlayer}{nodelayer}
		\node [style={ket_s}] (0) at (0, -0.25) {$\;x\;$};
		\node [style=none] (2) at (0, 0.5) {};
	\end{pgfonlayer}
	\begin{pgfonlayer}{edgelayer}
		\draw [style=x-wire] (0) to (2.center);
	\end{pgfonlayer}
\end{tikzpicture}
 \quad=\quad \sum_{i = 1}^k \enspace \input{diagrams/pos-elems-naive.tikz}\enspace.
\end{equation*}
We can get rid of the sum by replacing the $y_i\colon \C \to X$ by a map $y \colon \C^k \to X$ such that $y\ket{i} = y_i$ and write
\begin{equation*}
  \begin{tikzpicture}
	\begin{pgfonlayer}{nodelayer}
		\node [style={ket_s}] (0) at (0, -0.25) {$\;x\;$};
		\node [style=none] (2) at (0, 0.5) {};
	\end{pgfonlayer}
	\begin{pgfonlayer}{edgelayer}
		\draw [style=x-wire] (0) to (2.center);
	\end{pgfonlayer}
\end{tikzpicture}
 \quad=\quad \input{diagrams/pos-elems.tikz}\enspace.
\end{equation*}
Conversely, for any element of this form with arbitrary auxiliary Hilbert space $Z$, plugging in the completeness relation with respect to a self-conjugate basis on the connecting wire recovers $\dim Z$ elements $x_i$ and $x_i^*$.
Similar ideas show that complete positivity of operators between quantum sets can also be nicely expressed graphically without resorting to realignment and the double-wire calculus. See \cite{heunen_categories_2019} for more details.

\begin{lemma}[Theorem 7.18 in \cite{heunen_categories_2019}]\label{def:cp}
  Let $X, Y$ be quantum sets. A map $\varphi \colon X \to Y$ is completely positive if and only if there exists an object $Z$ and morphisms $g \colon X \otimes Y \to Z$, $h \colon Y \otimes X \to Z$, such that either of the following equivalent conditions holds.
  \begin{equation*}
    \input{diagrams/cp.tikz} \qquad\qquad\qquad \input{diagrams/cp2.tikz}\;.
  \end{equation*}
\end{lemma}

We will need a notion of subsets of quantum sets in order to define cuts, as well as source and sink sets for flows. A well-known approach defines subsets of a quantum set as the quotients of its underlying \Cstar-algebra. This definition is much too rigid for our purposes. Concretely, it can be shown that under this definition, $M_n$ has no non-trivial subsets. Instead we will work with the \emph{projective subsets} of \cite{spitzer_unification_2026}.

\begin{definition}
  Let $X$ be a quantum set. A \emph{projective subset} $Y \psubseteq X$ is equivalently given by
  \begin{enumerate}
    \item A projector $\mmap{Y} \colon X \to X$ satisfying
      \begin{equation*}
        \input{diagrams/right-module-cond.tikz}\enspace,
      \end{equation*}
    \item A projection $y \in X$,
    \item A right-module of $B(\C^N, \C)$ over the commutant of $X$, when viewing $X$ as a $*$-subalgebra of $B(\C^N)$. We denote the projector onto this module by $\pproj{Y}\colon B(\C^N, \C) \to B(\C^N, \C)$.
  \end{enumerate}
\end{definition}

If $X = M_n$, we have $X = B(\C^n)$ and $X' = \C\identity$. Further, by the Riesz representation theorem, we may identify $B(\C^n, \C)$ with $\C^n$. It follows that the projective subsets of $M_n$ are precisely given by projectors $\pproj{Y} \colon \C^n \to \C^n$. It can then be shown that in the double-wire calculus, we have
\begin{equation}\label{eq:pprojtommap}
  \input{diagrams/pproj-to-mmap.tikz}\enspace.
\end{equation}
For arbitrary quantum sets, $\mmap{Y}$ and $y$ are related through
\begin{equation*}
  \input{diagrams/mmap-to-elem.tikz} \qquad\quad\qquad \input{diagrams/elem-to-mmap.tikz}\enspace.
\end{equation*}
For classical sets, the definition recovers the usual notion of subsets. Concretely, the first condition requires that $\mmap{Y}$ project onto a subspace spanned by a subset of the standard basis vectors. Further, note that since $\mmap{Y}$ is a projector $\C^n \to \C^n$, we may also see it as defining a projective subset $Z \psubseteq M_n$ by letting $\pproj{Z} = \mmap{Y}$. We call subsets of this form \emph{classical subsets}. This idea can be extended to arbitrary quantum sets. Indeed, since $\mmap{Y}$ projects onto a coordinate subspace of $\C^n$, it projects onto a right-module over $\C^n$. Since for all \Cstar-algebras $X$ we have $X' \subseteq \C^n$, it follows that $\mmap{Y}$ is always of the form required in condition 3.

\begin{lemma}
  Let $X$ be a quantum set and $Y$ be a classical set. A projective subset $Z \psubseteq Y$ defines a projective subset $Z' \psubseteq X$ by setting $\pproj{Z'} = \mmap{Z}$. 
\end{lemma}

There is some ambiguity when defining the cardinality of a projective subset. Depending on the use case, it may be more natural to define it either as $\rk \mmap{Y}$ or $\rk \pproj{Y}$, see \cite[Remark 3.13]{spitzer_unification_2026}. Here we choose the former option.

\begin{definition}
  The \emph{cardinality} of a projective subset $Y \psubseteq X$ is given by $\abs{Y} = \rk \mmap{Y}$. Equivalently, its cardinality is given graphically by
  \begin{equation*}
    \abs{Y} \quad=\quad \input{diagrams/projsubset-card.tikz}\enspace. 
  \end{equation*}
\end{definition}

We illustrate these notions for the two basic examples of quantum sets.

\begin{example}\label{ex:projsubsetscn}
  If $X = \C^n$, an element $y = \sum_{i = 1}^n y_i \ket{i}$ is a projection, i.e.\ $y^2 = y = y^*$, if and only if $y_i^2 = y_i = y_i^*$ for all $i \in [n]$, since multiplication and involution in $\C^n$ act entrywise; that is, if and only if $y \in \{0, 1\}^n$. Projective subsets $Y \psubseteq \C^n$ are therefore in bijection with subsets $S \subseteq [n]$ via $y = \mathbf{1}_S \coloneqq \sum_{i \in S} \ket{i}$. By the relation between $y$ and $\mmap{Y}$ above, $\mmap{Y}$ is multiplication by $y$, so
  \begin{equation*}
    \mmap{Y} = \diag(y) = \sum_{i \in S} \ket{i}\bra{i}
  \end{equation*}
  is the projector onto $\operatorname{span}\{\ket{i} : i \in S\}$. The cardinality of $Y$ is thus
  \begin{equation*}
    \abs{Y} = \rk \diag(y) = \abs{S} = \sum_{i = 1}^n y_i,
  \end{equation*}
  the number of entries of $y$ equal to $1$. In particular, the full projective subset $\C^n \psubseteq \C^n$, given by $y = \mathbf{1}$ and $\mmap{Y} = \identity$, has cardinality $n$.
\end{example}

\begin{example}\label{ex:projsubsetsmn}
  Let $X = M_n = B(\C^n)$. A projection $y \in M_n$ is an orthogonal projector $P \colon \C^n \to \C^n$, and by the relation between $y$ and $\mmap{Y}$ above, $\mmap{Y}$ is left multiplication by $P$: $\mmap{Y}(x) = Px$ for all $x \in M_n$. We have $\mmap{Y} = P \otimes \identity$ and $\pproj{Y} = P$. The cardinality of $Y$ is thus
  \begin{equation*}
    \abs{Y} = \rk(P \otimes \identity) = n \rk P = n \rk \pproj{Y}.
  \end{equation*}
  In particular, the full projective subset $M_n \psubseteq M_n$, given by $P = \identity$, has cardinality $n^2 = \dim M_n$, while the projective subset given by a rank-one projector $P = \ket{\psi}\bra{\psi}$ has cardinality $n$, even though $\rk \pproj{Y} = 1$.
\end{example}

Finally, we recall the definition of quantum graphs. 

\begin{definition}[cf. Definition V.1 in \cite{musto_compositional_2018}]\label{def:qgraphsadj}
  A \emph{quantum graph} is a tuple $(X, G)$, where $X$ is a quantum set, called the \emph{vertex set}, and $G\colon X \to X$ is a linear map, called the \emph{adjacency operator}, satisfying
  \begin{equation*}
    \input{diagrams/adj-op-schur-idemp.tikz} \qquad\qquad\qquad \input{diagrams/adj-op-real.tikz}\enspace.
  \end{equation*}
\end{definition}

The conditions say that $G$ is a Schur idempotent, which classically corresponds to $G$ having entries in $\{0, 1\}$, and that $G$ is self-conjugate, or \emph{real}, which is classically implied by the first condition but is required for other quantum sets. 
Unless otherwise specified, we assume quantum graphs to be undirected and loopless. This corresponds to the conditions
\begin{equation*}
  \input{diagrams/adj-op-undirected.tikz} \qquad\quad\qquad \input{diagrams/adj-op-irreflexive.tikz}.
\end{equation*}

Quantum graphs may be equivalently defined as a vertex set $X$ together with a projective subset $E \psubseteq X \otimes X\op$ called the \emph{edge set}. This mirrors the classical case, where a subset $E \subseteq X \times X$ defines a graph with vertex set $X$, by assigning an edge from $u \in X$ to $v \in X$ whenever $(u, v) \in E$. To see that this definition is equivalent, note that $\pproj{E}$ is a projector and thus $\pproj{E}^R$ is a real Schur idempotent by \cref{prop:rlgnschurtocomp,prop:rlgnrealselfadj}.\footnote{Technically, we would also have to show that $\pproj{E}$ and $\pproj{E}^R$ have the correct types. Moreover, we have defined the realignment only for $X = M_n$. The correspondence holds for arbitrary quantum vertex sets. For more details refer to \cite{spitzer_unification_2026}.} Note that the image of $\pproj{E}$ is a subspace of $M_n$, which we call the \emph{quantum edge space} $\mathcal{S}_G$. It can be shown that $\mathcal{S}_G$ is closed under adjoints if and only if $G$ is undirected, and contains the identity if and only if $G$ has loops at every vertex. This recovers the operator system formalism of quantum graphs first introduced in \cite{duan_zeroerror_2013}.

\section{Quantum Cuts and Flows}\label{sec:cutsflows}
We can now start by introducing the fundamental definitions underlying our theory of cuts and flows in the quantum setting. Our definitions are motivated by the standard non-commutative dictionary.

\subsection{Capacity Maps}\label{sec:capacity-maps}
Common to both cuts and flows is the requirement for a notion of \emph{edge capacities}.
For undirected classical graphs, a capacity is a nonnegative number $c_{ij}$
assigned to an edge $\{i,j\}$. Extending $c$ by zero to all ordered pairs $(i,
j)$ then gives a symmetric matrix $C$ with non-negative entries satisfying $A
\bullet C=C$, where $A$ is the $\{0, 1\}$-valued adjacency matrix and $\bullet$
is the Schur product. Using the diagrammatic formalism, this definition has a
straightforward quantum analogue.

\begin{definition}
Let $(X, A)$ be a quantum graph. A \emph{capacity map} on $(X, A)$ is a
completely positive linear map $K\colon X\to X$ satisfying
\begin{equation*}
  \input{diagrams/capac-schur-cond.tikz}
  \qquad\qquad\qquad
  \input{diagrams/capac-sym-cond.tikz}
\end{equation*}
\end{definition}

We will sometimes call the first condition the \emph{support condition}, and say that $K$ is \emph{supported on $A$}. The condition that the capacities be non-negative is captured by the
requirement that $K$ be \emph{completely positive}. For classical graphs, that
is when the quantum set is $\C^n$, this condition is precisely equivalent to
$K$ being entrywise non-negative.

\begin{example}
  An important example of a capacity map is the adjacency operator $A_G$ of a
  quantum graph itself. We call this the \emph{unit capacity}. It corresponds
  to assigning a capacity of $1$ to each edge. It is clear that the adjacency
  operator is a valid capacity map: The support condition corresponds to
  the Schur idempotency condition for quantum graphs, and the adjacency
  operator is symmetric if the quantum graph is undirected. Finally, it is not
  hard to show that
  \begin{equation*}
    \input{diagrams/ag-edgeproj.tikz}
  \end{equation*}
  is a projector and hence $A_G$ is completely positive.
\end{example}

\subsection{Flows}

We can now move on to define flows on quantum graphs. We begin by recalling the
definition of flows on classical graphs. Recall that we can phrase a capacity
map on a classical graph as a symmetric matrix $(c_{ij})_{i, j \in V(G)}$
with non-negative entries such that $c_{ij} = 0$ whenever $i \not\sim j$.
A flow on a classical graph $G$
goes from a set $S \subseteq V(G)$ of \emph{source} vertices to a set $T
\subseteq V(G)$ of \emph{sink} vertices. The aim is to transport some commodity from the sources to the sinks while respecting the edge capacities; one is usually interested in the maximum possible flow from source to sink.

\begin{definition}[cf. Section 6.2 in \cite{diestel_graph_2017}]
  Let $G$ be a classical graph. Let $S, T \subseteq V(G)$ be disjoint sets and
  let $(c_{ij})_{i, j \in V(G)}$ be a matrix of capacities for $G$. An $S$-to-$T$
  \emph{flow} on $G$ is given by a function $f\colon V(G) \times V(G) \to \R$
  satisfying the following conditions.
  \begin{enumerate}[label=(CF\arabic*),leftmargin=5em]
    \item Orientation: $f(i, j) = -f(j, i)$.
    \item Capacity bound: $f(i, j) \leq c_{ij}$.
    \item Conservation: $\sum_k f(i, k) = 0$ for all $i \not\in S \cup T$.
  \end{enumerate}
\end{definition}

The first condition says that $f$ assigns flows to directed edges of the
graph such that opposite orientations get the same flow with opposite sign. The
second condition says that the flow on a single edge cannot exceed its capacity.
Finally, the third condition together with the first says that the flow has
no ``leaks'': the total flow into a vertex equals the flow out of it, unless 
the vertex is a source or a sink.
Note that the first and the second conditions above imply that a flow $f$ is supported on the graph: if $i$ and $j$ are non-adjacent, then $f(i, j) = 0$. We can then define the \emph{value} of the flow as the
net flow into the sink:
\begin{equation*}
  \val(f) = \sum_{t \in T} \sum_{k \not\in T} f(k, t).
\end{equation*}
Since flows have no leaks, the value may equivalently be defined as the net
flow out of the source.

\medskip

Now let $G$ be a quantum graph and let $K \colon V(G) \to V(G)$ be a capacity
map as defined in \cref{sec:capacity-maps}. To define flows of $G$, we first need to define quantum
analogues of the source and sink sets. The natural choice is to use the
projective subsets of \cite{spitzer_unification_2026}. Consider then projective subsets $S, T
\psubseteq V(G)$. Recall that they are in bijective correspondence with
projections $s, t \in V(G)$, and projectors $\mmap{S}, \mmap{T}\colon V(G) \to
V(G)$ that are also right-module homomorphisms, which are given by
\begin{equation*}
  \mmap{S}=L_s \qquad\quad \mmap{T}=L_t.
\end{equation*}
Disjointness of $S$ and $T$ then corresponds to $st = 0$ or equivalently
$\mmap{S} \circ \mmap{T} = 0$. In that case, $S$ and $T$ induce a third subset
$R \psubseteq V(G)$ given by
\begin{equation*}
  R \coloneqq V(G) \setminus (S \cup T) \qquad\quad \mmap{R} = \identity - \mmap{S} - \mmap{T}
\end{equation*}
such that $R, S, T$ partition $V(G)$.

The crucial observation to quantise the definition of flows is that similarly
to capacity maps, we can view a classical flow $f \colon V(G) \times V(G) \to
\R$ as a real-valued matrix $F$ satisfying additional conditions. The orientation condition corresponds to antisymmetry of
$F$. The capacity bound is classically equivalent to $K - F$ being entrywise
non-negative, which we translate again as complete positivity of the
corresponding operator. Finally, to generalise the conservation condition, we
note that the sum condition (CF3) precisely says that the $i$th row sum of $F$
is equal to 0 for all $i \in R$. Classically, the vector of row sums is given
by $F\mathbf{1}$ which we can restrict by projecting onto $R$. The all-ones
vector is precisely the unit of the classical set $\C^n$, so it is natural to
generalise this condition as $\mmap{R}Fu = 0$. We obtain the following
definition. For brevity, we give the conditions in text form, but it should be
clear how they are represented graphically.

\begin{definition}
Let $G$ be a quantum graph with adjacency operator $A_G$, and let $S, T
\psubseteq V(G)$ be disjoint. Let $K$ be a capacity map on $G$. A \emph{quantum
$S$-to-$T$ flow} on $(G, K)$ is a linear map $F\colon V(G) \to V(G)$ satisfying
the following conditions:
\begin{enumerate}[label=(F\arabic*),leftmargin=5em]
  \item Reality: $\overline{F}=F$.
  \item Antisymmetry: $F\transpose=-F$. Together with reality, this implies $F\adjoint=-F$.
  \item Capacity bound: $K-F$ is completely positive.
  \item Conservation: $\mmap{R} F u=0$.
\end{enumerate}
\end{definition}

The \emph{value} of a flow is its net divergence into the sink.
Classically, this corresponds to the sum of the $i$th column sums for $i \in T$.
In terms of matrices this value is given by $\langle \mathbf{1}, FT\mathbf{1}
\rangle$. We again generalise this in terms of the unit which gives
\begin{equation*}
  \val(F)
  \coloneqq
  (u\adjoint F \mmap{T} u)(1)
  \quad=\quad \input{diagrams/flow-value.tikz}\enspace.
\end{equation*}
Crucially, $\val(F)$ is a real number. To see this, we note that by
\cite[Lemma 3.3]{spitzer_unification_2026}, we have $t = \mmap{T}u$, which satisfies
$t^* = t$. Since the \Cstar-involution corresponds to the conjugation functor
in $\FdHilb$, we have
\begin{equation*}
  \overline{\mmap{T}u} = t^* = t = \mmap{T}u.
\end{equation*}
It now follows that
\begin{equation*}
  \overline{u\adjoint F \mmap{T}u} = u\adjoint \overline{F} \overline{\mmap{T}u} = u\adjoint F \mmap{T} u,
\end{equation*}
and thus $\val(F)$ must be real.

Using the properties of the flow, one can also show that its value is equal to the net
divergence out of the source. Indeed, because $F\adjoint=-F$, one has
$u\adjoint F u = 0$. Using conservation,
\begin{equation*}
  0
  =
  u\adjoint F u
  =
  u\adjoint \mmap{S}F u
  +
  u\adjoint \mmap{T}F u.
\end{equation*}
Moreover, since $\val(F)$ is a real scalar, we can take its adjoint without changing its
value. $F\adjoint=-F$ and $\mmap{T}\adjoint=\mmap{T}$ then give
\begin{equation*}
  u\adjoint F\mmap{T} u
  =
  (u\adjoint F\mmap{T} u)\adjoint
  =
  - u\adjoint \mmap{T}F u.
\end{equation*}
Combining the two equations, we get $u\adjoint \mmap{S} F u = u\adjoint F
\mmap{T} u$, so the source-outflow and sink-inflow expressions agree.

Finally, we can consider the \emph{maximum} flow from a source set $S$ to a sink set $T$.

\begin{definition}
Let $G$ be a quantum graph, $K$ a capacity on $G$, and $S, T \psubseteq V(G)$
disjoint. The \emph{maximum $S$-to-$T$ flow} value is
\begin{equation*}
  \MaxFlow(G,K;S,T)
  \coloneqq
  \max\bigl\{
    \val(F)
    \mid
    F\text{ is an }S\text{-to-}T\text{ flow on }(G,K)
  \bigr\}.
\end{equation*}
If $K$ is the unit capacity, we will generally omit it from the notation.
\end{definition}

\subsection{Cuts}

It remains to define cuts of quantum graphs. Classically, any binary partition
of the vertex set uniquely defines a cut by considering the edges between the
two blocks. If we are moreover given a capacity on our graph, then the cost of
this cut is given by the sum of the capacities of the edges. We can quantise
this definition again in terms of projective subsets. We first recall the
definition of partitions for quantum sets.

\begin{definition}[cf. Definition 4.5 in \cite{spitzer_unification_2026}]
  Let $X$ be a quantum set. A \emph{partition} of $X$ into $k$ blocks is given by a collection of projective subsets $P_1, \dots, P_k \psubseteq X$ such that $\sum_{i = 1}^k \mmap{P_i} =\I $.
\end{definition}

This leads to the following quantum generalisation of cuts.

\begin{definition}
  Let $G$ be a quantum graph. A \emph{cut} is a partition of $V(G)$ into two subsets $E, F \psubseteq V(G)$, that is
\begin{equation*}
  \mmap{E}+\mmap{F}=\I,\qquad
  \mmap{E} \mmap{F} = \mmap{F}\mmap{E} = 0.
\end{equation*}
If $S, T \psubseteq V(G)$ then $(E, F)$ is an \emph{$S$-$T$-cut} if $S \psubseteq E$ and $T \psubseteq F$, that is
\begin{equation*}
  \mmap{S}\mmap{E}=\mmap{E}\mmap{S} = \mmap{S},\qquad
  \mmap{T}\mmap{F}=\mmap{F}\mmap{T} = \mmap{T}.
\end{equation*}
If $K$ is a capacity on $G$, then the \emph{cost} of the cut is
\begin{equation*}
  \cut_K(E,F)
  \coloneqq
  (u\adjoint \mmap{E} K \mmap{F} u)(1).
\end{equation*}
\end{definition}

In the classical case, the cost definition reduces precisely to the sum of edge
capacities across the cut. Indeed, we have seen that a capacity map on $\C^n$
is just an entrywise non-negative symmetric matrix, and projective subsets of
classical sets are given by diagonal projectors on $\C^n$. The classical cut
value $\cut^{\mathrm{cl}}_{K}$ is thus given by summing over the entries of the
top-left (equivalently bottom-right) quadrant of the capacity matrix,
$\cut^{\mathrm{cl}}_K = \langle \mathbf{1}, \mmap{E}K\mmap{F}\mathbf{1}\rangle$.
Since $\mathbf{1} = u(1)$, this is precisely equal to $\cut_K(E, F)$.

It follows again from the properties of projective subsets that the cost of a
cut is a nonnegative real number. Recall that $e = \mmap{E}u$ and $f =
\mmap{F}u$ are projections, so $e = e^* = e^2$ and $f = f^*f$ is a positive
element of $V(G)$. Since $K$ is completely positive, it sends positive elements
to positive elements. By the graphical characterisation of positive elements
recalled in \cref{sec:prelims}, we can thus write
\begin{equation*}
  \input{diagrams/kf-positive.tikz}
\end{equation*}
and we get
\begin{equation*}
  \cut_K(E, F) \quad=\quad \input{diagrams/cut-cost-geq0.tikz}.
\end{equation*}
The corresponding \emph{minimum $S$-$T$-cut value} is denoted by
\begin{equation*}
  \MinCut(G,K;S,T)
  \coloneqq
  \inf
  \bigl\{
    \cut_K(E,F)
    \mid
    E,F \psubseteq V(G),\ \mmap{E} + \mmap{F} = \I,\ \mmap{S}\mmap{E} = \mmap{S},\ \mmap{T}\mmap{F}=\mmap{T}
  \bigr\}.
\end{equation*}

\begin{remark}
We remark that the cut is of course already uniquely defined by a single
subset $E \psubseteq V(G)$ since $\mmap{F} = \I - \mmap{E}$, but we will always state both
subsets to simplify the formulas.
\end{remark}

We will also consider the following generalisation of cuts on classical graphs.
An ordinary $S$-$T$-cut is uniquely defined by a subset $E$ of the vertex set
such that $S \subseteq E$ and $T$ and $E$ are disjoint. Instead of vertices
belonging to $E$ or not, we can generalise this definition by assigning to each
vertex a \emph{weight} with which it contributes to $E$. We require vertices in
$S$ to contribute fully, while vertices of $T$ must not contribute at all. The
cost of the generalised cut is computed by multiplying the capacities of edges
leaving a vertex by its weight. If both endpoints have non-zero weight their difference matters.
This definition is also known in the literature
as the Lovász extension of the cut function \cite[Section 6.2]{bach2011learning}
and has applications in graph neural networks and graph signal processing \cite{sardellitti2017graph}.

\begin{definition}(Generalised Cuts)\label{def:gencuts}
  Let $G$ be a classical graph and let $S, T \subseteq V(G)$, $S \cap T = \emptyset$. Let $C = (c_{ij})_{i,
  j \in V(G)}$ be a (symmetric) capacity matrix for $G$. A \emph{generalised $S-T$ cut}
  is given by a weight function $w\colon V(G) \to \R$ such that
  \begin{enumerate}
    \item $w(s) = 1$ for all $s \in S$,
    \item $w(t) = 0$ for all $t \in T$.
  \end{enumerate}
  Its \emph{cost} is given by
  \begin{equation*}
    \operatorname{gen-cut}_C(w) \coloneqq \sum_{i < j \in V(G)} c_{ij} \cdot \abs{w(i) - w(j)}.
  \end{equation*}
\end{definition}

Clearly, an ordinary cut on a classical graph is also a generalised cut with
$w(i) = [i \in E]$, in which case its cost is given by
\begin{align*}
  \operatorname{gen-cut}_C(w) &= \sum_{i < j \in V(G)} c_{ij} \cdot \abs{[i \in E] - [j \in E]}\\
                              &= \sum_{i, j \in V(G)} c_{ij} \cdot [i \in E \land j \not\in E]\\
                              &= \langle \mathbf{1}, EC(\I - E) \mathbf{1} \rangle\\
                              &= \cut_C(E, \I - E).
\end{align*}

A natural question is whether we can extend this definition to quantum graphs.
We quantised ordinary cuts by passing from subsets to projective subsets. These
are in bijective correspondence with projections $f \in V(G)$ through the
equation $\mmap{F} = L_f$. It was argued in \cite{spitzer_unification_2026} that $f$ is
the analogue of the \emph{characteristic vector} of the subset. Similarly to
the idea of a characteristic vector, a function $w\colon [n] \to \R$ can be
represented by a real vector $w' \in \C^n$, where $w'_i = w(i)$. These real vectors
are precisely the self-involutive elements of $\C^n$. Following this line of
thinking, we may quantise functions $V(G) \to \R$ to be the self-involutive
elements of $V(G)$. Through left multiplication, we can also lift these real
elements to maps $V(G) \to V(G)$. The maps we obtain this way are precisely the
self-adjoint right-module homomorphisms.

\begin{proposition}
  Let $X$ be a quantum set. There is a bijective correspondence between
  elements $w \in X$ satisfying $w^* = w$ and self-adjoint right-module
  homomorphisms $\mmap{W} \colon X \to X$.
\end{proposition}

The proof is the same as for Lemma 3.3 in \cite{spitzer_unification_2026}.

\begin{definition}
  Let $X$ be a quantum set. A \emph{real-valued function} $F\colon X
  \rightsquigarrow \R$ on $X$ is equivalently given by a self-involutive
  element $f \in X$ or a self-adjoint right-module homomorphism $\mmap{F}\colon
  X \to X$.
\end{definition}

\begin{remark}
  In the spirit of Musto, Reutter, and Verdon's quantum functions \cite[Definition III.11]{musto_compositional_2018},
  one could easily generalise this definition to \emph{quantum real-valued
  functions} by graphically adding a Hilbert space wire through $f$ and
  $\mmap{F}$, such that $\mmap{F}$ becomes a map $\Hilb \otimes X \to X \otimes
  \Hilb$ and $f$ a function $\Hilb \to X \otimes \Hilb$. We will not need this
  generalisation here, but it might be an interesting object of future study.
\end{remark}

With this definition in hand, we can quantise the concept of the weight
function. The conditions on $S$ and $T$ can be unified in the single condition
\begin{equation*}
  (\mmap{S} + \mmap{T})\mmap{W}(\mmap{S} + \mmap{T}) = \mmap{S},
\end{equation*}
which says that when restricted to $S \cup T$, $\pproj{W}$ is just the indicator function of $\pproj{S}$.
This is already enough to define the
generalised cut itself, but it remains to generalise the definition of its
cost. The goal is to again phrase it in terms of a sum over matrix elements,
which we can then generalise in terms of operators and inner products. Indeed,
this turns out to be possible.

\begin{lemma}\label{lem:cutcostmin}
  Let $G$ be a classical graph and let $K$ be a capacity on $G$. Let $w\colon V(G) \to \R$ be a weight function defining a generalised cut. Then we have
  \begin{equation*}
  \operatorname{gen-cut}_K(w) = \min_{\substack{B = \overline{B} = B\transpose\\ \breve{W} - B \mathrm{~entrywise~nonneg.}}} \langle \mathbf{1}, ((\breve{W} - B) \bullet K) \mathbf{1} \rangle
  \end{equation*}
  where $\breve{W} \in \R^{V(G) \times V(G)}$ with $\breve{W}_{ij} = w(i)$.
  \begin{proof}
    Consider a pair of matrix elements $\breve{W}_{ij}$ and $\breve{W}_{ji}$. We can always set one of the two values to $0$ by subtracting the symmetric operator
    \begin{equation*}
       B^{(i,j)} = \breve{w}^{(i, j)} \ket{i}\bra{j} + \breve{w}^{(i, j)} \ket{j}\bra{i},
    \end{equation*}
    where  $\breve{w}^{(i, j)} = \min \{\breve{W}_{ij}, \breve{W}_{ji} \}$. In this case, the value of the remaining entry becomes
    \begin{equation*}
      \abs{\breve{W}_{ij} - \breve{W}_{ji}} = \abs{w(i) - w(j)}.
    \end{equation*}
    We now let
    \begin{equation*}
      B = \sum_{i < j} B^{(i, j)} + \sum_i \breve{W}_{ii} \ket{i}\bra{i}.
    \end{equation*}
    Evaluating the objective function then yields
    \begin{equation*}
      \langle \mathbf{1}, ((\breve{W} - B) \bullet K) \mathbf{1} \rangle
      =
      \sum_{i < j} K_{ij} \abs{w(i) - w(j)}
    \end{equation*}
    as desired.

    It is also clear that we cannot do better: The objective is equal to a sum
    of nonnegative elements, weighted by the entries of $\breve{W}' = \breve{W}
    - B$. The diagonal elements $\breve{W}'_{ii}$ are zero by construction. It
    follows that to reduce the objective value, we have to reduce at least one
    coefficient $\breve{W}'_{ij}$ for $i \neq j$. Without loss of generality,
    we have $\breve{W}'_{ij} > 0$, which implies $\breve{W}'_{ji} = 0$ by
    construction. But the only way to decrease $\breve{W}'_{ij}$ is to increase
    $B_{ij}$, in which case we get $\breve{W}'_{ji} = \breve{W}_{ji} - B_{ji} =
    \breve{W}_{ji} - B_{ij} < 0$, leaving the feasible region.
  \end{proof}
\end{lemma}

If we let $w' \in \C^n$ be the vector representation of $w$, that is $w'_i =
w(i)$, we can write $\breve{W} = w'\mathbf{1}\adjoint$. We have seen that the quantum
analogue of this vector representation is precisely the representation of a
real-valued function $X \rightsquigarrow \R$ as a self-involutive element $x
\in X$. We have also seen that the correct quantisation of an entrywise
nonnegative matrix is a completely positive operator. We thus obtain the
following definition of generalised cuts on quantum graphs as a straightforward
quantisation of \cref{def:gencuts} and \cref{lem:cutcostmin}.

\begin{definition}
  Let $G$ be a quantum graph, and let $S, T \psubseteq V(G)$ be disjoint.
  Moreover, let $K$ be a capacity on $G$. A \emph{generalised $S-T$ cut} is
  given by a weight function $W\colon V(G) \rightsquigarrow \R$ such that
  \begin{equation*}
    (\mmap{S} + \mmap{T})\mmap{W}(\mmap{S} + \mmap{T}) = \mmap{S},
  \end{equation*}
  Its \emph{cost} is given by
  \begin{equation*}
    \operatorname{gen-cut}_K(W) \coloneqq \min_{\substack{ B = \overline{B} = B\transpose\\ \breve{W} - B \mathrm{~CP}}} [u\adjoint \circ ((\breve{W} - B) \bullet K) \circ u](1)
  \end{equation*}
  where $\breve{W} = w \circ u\adjoint$.
\end{definition}

In the classical case, we can solve for $B$ and hence get rid of the
minimisation, but we will see later that this is generally impossible for
arbitrary quantum sets. It remains to show that $\operatorname{gen-cut}_K(W)$
is well-defined by showing that the minimisation term is always real. Indeed, we have
\begin{equation*}
  \input{diagrams/gen-cut-real.tikz}
\end{equation*}
since $\breve{W} - B$ is completely positive in the feasible region and every
completely positive map is real. We can now define the corresponding
\emph{generalised minimum $S$-$T$ cut} cost as
\begin{equation*}
  \GenMinCut(G, K; S, T) \coloneqq \inf \{\operatorname{gen-cut}_K(W) \mid W\colon V(G) \rightsquigarrow \R,  (\mmap{S} + \mmap{T})\mmap{W}(\mmap{S} + \mmap{T}) = \mmap{S} \}.
\end{equation*}

\pagebreak
\section{The Case of \texorpdfstring{$M_n$}{Mn}}\label{sec:mn}\label{sec:mncase}
To dive deeper into the theory of flows and cuts, we now restrict ourselves to
quantum graphs whose vertex set is $M_n$. We uncover a rich theory that
beautifully quantises the classical landscape. We find that while classical
flows are described by linear programs, which famously yields the max-flow
min-cut theorem, flows on quantum graphs are described by semidefinite
programs, which we will use to derive a quantum version of the max-flow min-cut
theorem. We note that our restriction to $M_n$ is merely for clarity of our
presentation: Stating the SDPs for arbitrary quantum sets---direct sums of
matrix algebras---introduces cross-terms that one has to tediously keep track
of. We prove the general case of the theorem in Appendix \ref{app:cutflowgeneral}.

\subsection{The Choi Picture}

Since the realignment is invertible, we can state equivalent characterisations
of capacity maps, flows, and cuts by realigning the corresponding conditions.
We call this the \emph{Choi picture}. We have seen that our choice of
realignment corresponds precisely to the map
\begin{equation*}
  \Phi \mapsto \frac{1}{n} \sum_{i, j} \ket{i}\bra{j} \otimes \Phi(\ket{i}\bra{j}).
\end{equation*}
It thus follows by Choi's theorem that if a map was required to be
completely positive, it is required to be positive semidefinite in the Choi
picture. This yields the equivalent conditions for the complete positivity of
the capacity map and, by linearity of the realignment, the capacity bound on
flows.

The next important observation is that the transpose commutes with the realignment.
This implies that the symmetry condition on capacity maps and the antisymmetry
condition on flows stay invariant when moving to the Choi picture.

\begin{lemma}\label{lem:trposerealigncommute}
  For any linear map $F \colon M_n \to M_n$ it holds that $(F^R)\transpose = (F\transpose)^R$.
  \begin{proof}
    Follows from isomorphism of diagrams:
    \begin{equation*}
      \input{diagrams/realignment-transpose-comm.tikz} \qedhere
    \end{equation*}
  \end{proof}
\end{lemma}

We have also seen that the realignment acts as a kind of ``Fourier transform''
between the Schur product and ordinary composition of linear maps. The support
condition of capacity maps is thus converted into a range condition:
Consider $\Pi_A=A^R$. Since $A$ is a Schur idempotent and real, $\Pi_A$ is a
projector and we have
\begin{equation*}
  A \bullet K=K
  \quad\Longleftrightarrow\quad
  \Pi_A K^R = K^R.
\end{equation*}
Since $K$ is completely positive, $K^R \succeq 0$ and consequently $K^R$ is self-adjoint.
The condition is thus equivalent to
\begin{equation*}
  K^R=\Pi_A K^R \Pi_A.
\end{equation*}
Thus a capacity is equivalently given by a positive operator whose support lies
inside the quantum edge space $\mathcal S_A$. Similarly, reality of the flow map
becomes self-adjointness of its realignment.
We derive the counterpart of the conservation condition for flows graphically. We have
\begin{equation*}
  0 \quad\stackrel{!}{=}\quad \input{diagrams/conservation-realigned.tikz}\enspace,
\end{equation*}
where we denoted $\pproj{R} \coloneqq \I - \pproj{S} - \pproj{T}$. 
We can then bend down the rightmost wire to obtain the equivalent condition $\pproj{R}\Tr_1 F^R = 0$.
Similarly, we can give the value of a flow $F$ in terms of $F^R$ as
\begin{equation*}
  \input{diagrams/flow-value-realigned.tikz}\enspace.
\end{equation*}

We summarise the correspondence between the original conditions and the
equivalent conditions in the Choi picture in the following table.
To simplify the normalisation, we let $Z \coloneqq n^{2}F^R$ and $Q \coloneqq n^{2}K^R$.

\begin{center}
\begin{tabular}{rl}
\toprule
\textbf{Map condition} & \textbf{Choi/realignment condition} \\
\midrule
Symmetry of \(K\) & \(Q=Q\transpose\) \\
Support on the graph & \(Q=\Pi_A Q \Pi_A\) \\
\(K\) completely positive & \(Q\) positive semidefinite  \\
\midrule
Reality of \(F\) & \(Z=Z\adjoint\) \\
Antisymmetry \(F\transpose=-F\) &
\(Z\transpose=-Z\) \\
Capacity bound \(K-F\) CP &
\(Q-Z\succeq 0\). \(Q+Z\succeq0\) follows from $Q\transpose = Q$, $Z\transpose = -Z$. \\
Conservation & \((I-\pproj{S}-\pproj{T})\Tr_1 Z=0\) \\
Value & \(\val(F)=\Tr\bigl[(\I_n\otimes \pproj{S})Z\bigr] = \lrangle{\I_n \otimes \pproj{S}, Z}\)  \\
\bottomrule
\end{tabular}
\end{center}

\subsection{Semidefinite Programming Duality}

The conditions obtained through realignment are all either linear equality constraints or positivity constrains. This suggests to phrase the maximum-flow problem as an SDP. With the above notation, the maximum-flow problem on $M_n$ is

\begin{equation}\label{eq:P}
\begin{array}{ll}
\text{maximize} &
\lrangle{\I_n \otimes \pproj{S}, Z} \\
\text{subject to} &
  Z=Z\adjoint, \\
& Z\transpose=-Z, \\
& \pproj{R}\Tr_1 Z=0, \\
& Q-Z\succeq 0
\end{array}
\tag{P}
\end{equation}
where we again denoted $\pproj{R} \coloneqq \I - \pproj{S} - \pproj{T}$.

This situation nicely mirrors the classical case, where the maximum-flow
problem can be phrased as a linear program. One can compute the dual of this
linear program, and find that it corresponds precisely to a minimum cut of the
graph. This yields the celebrated max-flow min-cut theorem
\cite{dantzig2003max}.

\begin{theorem}[Max-Flow Min-Cut]
  Let $G$ be a classical graph, $K$ a capacity on $G$ and $S, T \subseteq V(G)$ be disjoint. Then $\MaxFlow(G, K; S, T) = \MinCut(G, K; S, T)$.
\end{theorem}

The natural next step is to compute the dual of the above-mentioned SDP for the
maximum flow problem and verify whether it corresponds to some notion of
minimum cut on quantum graphs.

\begin{proposition}
  The dual of the SDP \eqref{eq:P} is given by
  \begin{equation}\label{eq:D}
    \begin{array}{ll}
      \text{minimize} &
      \lrangle{\I_n \otimes (\pproj{S} - \tilde\Lambda_2) - \tilde\Lambda_1, Q} \\
      \text{subject to} & \tilde\Lambda_1\adjoint = \tilde\Lambda_1,
                          \tilde\Lambda_1\transpose = \tilde\Lambda_1,\\
                        & \tilde\Lambda_2\adjoint = \tilde\Lambda_2,\\
                        & \pproj{R}^\bot\tilde\Lambda_2\pproj{R}^\bot = 0,\\
                        & \I_n \otimes (\pproj{S} - \tilde\Lambda_2) - \tilde\Lambda_1  \succeq 0
    \end{array}
    \tag{D}
  \end{equation}
  \begin{proof}

    We assign the dual variables $\Lambda_1, \Lambda_2, \Lambda_3$ to the last three constraints in \eqref{eq:P} from top to bottom.
The Lagrangian is then given by
\begin{equation*}
  \mathcal{L}(\Lambda_1, \Lambda_2, \Lambda_3) = \lrangle{\I_n \otimes \pproj{S}, Z} -  \lrangle{Z\transpose + Z, \Lambda_1} - \operatorname{Re}\; \lrangle{\pproj{R} \Tr_1 Z, \Lambda_2} + \lrangle{Q - Z, \Lambda_3},
\end{equation*}
noting that $\pproj{R}\Tr_1 Z$ is not necessarily self-adjoint.
The value of the primal is now upper-bounded by
\begin{align*}
  &\min_{\substack{\Lambda_1\adjoint = \Lambda_1,\Lambda_2\\\Lambda_3 \succeq 0}}
  \max_{Z} ~ \mathcal{L}(\Lambda_1, \Lambda_2, \Lambda_3)\\
  =
  &\min_{\substack{\Lambda_1\adjoint = \Lambda_1, \Lambda_2 \\\Lambda_3 \succeq 0}}
  \max_{Z} ~ \lrangle{\Lambda_3, Q} + \lrangle{\I_n \otimes \pproj{S} - (\Lambda_1 + \Lambda_1\transpose) - \I_n \otimes (\pproj{R}\Lambda_2 + \Lambda_2\adjoint \pproj{R}) / 2 - \Lambda_3, Z} \\
  =
  &\min_{
    \substack{
      \tilde{\Lambda}_1\adjoint = \tilde{\Lambda}_1,
      \tilde\Lambda_2\adjoint=\tilde\Lambda_2,\\
      \pproj{R}^\bot\tilde\Lambda_2\pproj{R}^\bot = 0,\\
      \tilde{\Lambda}_1\transpose = \tilde{\Lambda}_1,\\
      \I_n \otimes (\pproj{S} - \tilde\Lambda_2) - \tilde{\Lambda}_1\succeq 0
    }
    }
    \lrangle{\I_n \otimes (\pproj{S} - \tilde\Lambda_2) - \tilde{\Lambda}_1, Q}.
\end{align*}
Here we wrote $\tilde\Lambda_1 = \Lambda_1 + \Lambda_1\transpose$ and $\tilde\Lambda_2 = (\pproj{R}\Lambda_2 + \Lambda_2\adjoint \pproj{R})/2$. This works since maps of these forms are precisely the symmetric maps and the maps satisfying $\pproj{R}^\bot X \pproj{R}^\bot = 0$, respectively. This yields the desired dual.
  \end{proof}
\end{proposition}

As is common with convex optimisation problems, we can show that \eqref{eq:D} satisfies strong duality.

\begin{proposition}
  The duality gap between the SDPs \eqref{eq:D} and \eqref{eq:P} is $0$.
  \begin{proof}
    By \emph{Slater's condition}, the duality gap is $0$ if the dual has a strictly
    feasible point and the primal is feasible and bounded below (see, e.g., \cite[Sec.~5.2.3]{boyd2004convex} or \cite[Theorem 1.18]{watrous2018theory}). Taking
    $\tilde\Lambda_2 = 0$ and $\Lambda_1 = -c\I$ with $c > \norm{\I \otimes S}$, the
    only inequality constraint
    \begin{equation*}
      \I \otimes (\pproj{S} - \tilde\Lambda_2) - \tilde{\Lambda}_1 = \I \otimes \pproj{S} + c\I \succ 0
    \end{equation*}
    is satisfied strictly. Taking $Z=0$ in \eqref{eq:P} shows that the primal is feasible
    and non-negative, finishing the proof.
  \end{proof}
\end{proposition}

We can further simplify the dual SDP. Note that since $\pproj{R}^\bot = \pproj{S} + \pproj{T}$ and $\pproj{S}\pproj{T} = \pproj{T}\pproj{S} = 0$, the maps of the form $\pproj{S} + \tilde\Lambda_2$ satisfy
\begin{equation*}
  (\pproj{S} + \pproj{T}) (\pproj{S} + \tilde\Lambda_2) (\pproj{S} + \pproj{T}) = \pproj{S},
\end{equation*}
and conversely every self-adjoint map satisfying this equation is of the above form. We can thus state the following equivalent SDP.

  \begin{equation}\label{eq:DD}
    \begin{array}{ll}
      \text{minimize} &
      \lrangle{\I \otimes \pproj{W} - \Lambda_1, Q} \\
      \text{subject to} & \Lambda_1\adjoint = \Lambda_1,\\
                        & \Lambda_1\transpose = \Lambda_1,\\
                        & \pproj{W}\adjoint = \pproj{W}, \\
                        & (\pproj{S} + \pproj{T})\pproj{W}(\pproj{S} + \pproj{T}) = \pproj{S},\\
                        & \I \otimes \pproj{W} - \Lambda_1 \succeq 0
    \end{array}
    \tag{D'}
  \end{equation}

  The choice of notation in $\pproj{W}$ is deliberate. Observe that $\pproj{W}$
  defines precisely a real-valued function $W\colon V(G) \rightsquigarrow \R$.
  Moreover, the condition
  $(\pproj{S} + \pproj{T})\pproj{W}(\pproj{S} + \pproj{T}) = \pproj{S}$
  is equivalent to 
  $(\mmap{S} + \mmap{T})\mmap{W}(\mmap{S} + \mmap{T}) = \mmap{S}$ by \cref{eq:pprojtommap}.

  Now let us express the objective graphically and
  rewrite it in terms of $K = n^{-2}Q^\bR$:
  \allowdisplaybreaks
  \begin{equation}\label{eq:dualcutinterp}
  \begin{aligned}
    \input{diagrams/dual-sdp-realigned.tikz} &\quad\enspace=\enspace\quad \input{diagrams/dual-sdp-realigned-2.tikz}\\
                                 &\quad\enspace=\enspace\quad \input{diagrams/dual-sdp-realigned-3.tikz}\\
                                 &\quad\enspace=\enspace\quad \input{diagrams/dual-sdp-realigned-4.tikz}
  \end{aligned}
  \end{equation}
  By \cref{lem:trposerealigncommute}, $\Lambda_1^\bR$ is symmetric if and only if $\Lambda_1$ is. Moreover, a quick diagrammatic derivation shows that $(wu\adjoint)^R = \I \otimes \pproj{W}$, so $wu\adjoint - \Lambda_1^\bR$ is completely positive if and only if $\I \otimes \pproj{W} - \Lambda_1 \succeq 0$. It follows that minimising over symmetric $\Lambda_1$ satisfying the PSD condition yields precisely
  \begin{equation*}
    \min_{\substack{\Lambda_1 = \Lambda_1\transpose\\ wu\adjoint - \Lambda_1^\bR \text{ CP}}} [u\adjoint \circ (K \bullet (wu\adjoint - \Lambda_1^\bR)) \circ u](1) = \operatorname{gen-cut}_K(W).
  \end{equation*}
  The full optimisation minimises over valid weight functions and thus 
  has optimum $\GenMinCut(G, K; S, T)$.

  We see that just as in the purely classical case, there exists a strong
  duality between maximum flows and minimum cuts. The main difference is that
  the duality is not with ordinary cuts but with generalised cuts. The
  attentive reader might notice that this seems to contradict our statement at
  the beginning of the section that the derivations go through for
  arbitrary quantum sets and not just $M_n$, since that would subsume the
  classical case. This tension is resolved by noting that the minimum
  generalised cut of a classical graph is equal to the ordinary minimum cut.
  This follows from a convexity argument, noting that we can optimize the
  classical weight function entrywise.

  \begin{proposition}\label{prop:clmincutisgenmincut}
    Let $G$ be a classical graph. Then $\MinCut(G, K; S, T) = \GenMinCut(G, K; S, T)$ for all choices of $K, S, T$.
  \end{proposition}

  In some sense, the duality with the generalised cut is even more natural than the
  classical duality with ordinary cuts: In the classical derivation, only the
  optimum of the dual linear program may be interpreted as a cut. On the other
  hand, our dual SDP can be interpreted naturally as a minimisation over all
  valid generalised cuts.

\subsection{On Ordinary Cuts}\label{ssec:ordcuts}

Since we have shown a duality between the maximum flow and the minimum generalised cut, a natural question is how ordinary cuts fit into this picture. We can show that, as one would expect, $\MinCut(G, K; S, T)$ provides an upper bound for $\GenMinCut(G, K; S, T)$. 

Recall that an ordinary cut is uniquely determined by the subset $E \psubseteq V(G)$ containing $S$, which in turn can be seen as a zero-one valued version of the weight function $W$ of a generalised cut. Indeed, since $e$ is a projection, it is in particular a self-involutive element of $V(G)$ and thus a valid weight function. 
Consider now the matrix $\breve{E} = e \circ u\adjoint$ used to define the cost of a generalised cut. In the classical case, $\breve{E}$ is a matrix whose $i$th column is the all-ones vector if $i \in E$ and the zero vector if $i \not\in E$. The largest symmetric matrix that we may subtract from $\breve{E}$ while keeping it entrywise non-negative is precisely $ee\adjoint$, the all-ones matrix supported on $E$. If we interpret $\breve{E}$ as a matrix indexed by vertices of $G$, then subtracting $ee\adjoint$ precisely sets those entries to zero that belong to edges between vertices in $E$. Consequently, taking the Schur product with $K$ and computing the sum of entries yields precisely the value of the cut.
This argument stays valid in the quantum case: By reality, $ee\adjoint$ is symmetric. Moreover, it holds that
\begin{equation*}
  \input{diagrams/ord-cut-cp.tikz}\enspace,
\end{equation*}
which is a projector, so $\breve{E} - ee\adjoint$ is completely positive.
Computing the analogue of the sum of entries of the Schur product with the capacity map then yields precisely
\begin{equation*}
  \input{diagrams/ord-cut-val.tikz} \quad=\quad \cut_K(E, F).
\end{equation*}
Now by \cref{eq:dualcutinterp}, the above is in turn equal to the objective of (\ref{eq:DD}) for $\pproj{W} = \pproj{E}$ and $\Lambda_1 = (ee\adjoint)^R$. Moreover, we have seen that $(eu\adjoint)^R = I \otimes \pproj{E}$, so these choices are feasible. It follows that $\MinCut(G, K; S, T)$ is at least as large as the SDP value and hence at least as large as $\GenMinCut(G, K; S, T)$.

\smallskip

We can now state the following quantum version of the max-flow min-cut theorem.

\thmmfmc
\begin{proof}
  We have shown the case of $V(G) = M_n$. The case of arbitrary quantum sets is proved in Appendix \ref{app:cutflowgeneral}.
\end{proof}

One might ask whether \cref{prop:clmincutisgenmincut} could be generalised to
quantum graphs, turning the last inequality into an equality. This is not
the case. In fact, there is already a simple separation between maximum flow
and projection-valued minimum cut in $M_2$.
Indeed, let
\begin{equation*}
  \pproj{S}=\ket0\bra0,
  \qquad
  \pproj{T}=\ket1\bra1
\end{equation*}
Then $\mmap{S} + \mmap{T}=\I$, so the conservation projection $\mmap{R}=\I-\mmap{S}-\mmap{T}$ is zero and
there is no interior conservation constraint. Consider the quantum graph $G$ whose Choi support is
\begin{equation*}
  \Pi_G
  =
  \lspan \{\ket{0}\bra{1}, \ket{1}\bra{0} \}
  \cong
  \lspan\{\ket{01}, \ket{10}\}
  \subseteq \C^2\otimes\C^2 .
\end{equation*}
A valid capacity on this graph is given by
\begin{equation*}
  Q(c,r)
  =
  c\,\ket{01}\bra{01}
  +c\,\ket{10}\bra{10}
  +r\,\ket{01}\bra{10}
  +r\,\ket{10}\bra{01} = \begin{bmatrix}
    c&r\\
    r&c
  \end{bmatrix},
  \qquad
  c\ge |r|.
\end{equation*}
The condition $c\ge |r|$ is exactly $Q(c,r)\succeq0$, i.e. complete
positivity of the capacity map.

The flow variable $Z$ is Hermitian, and antisymmetric under the categorical transpose $Z\transpose$.
In the standard basis, this involution is given by $(-)\transpose \colon X\mapsto \Sigma X^t \Sigma$.
It follows that the matrix representation of $Z$ in this basis must be of the form
\begin{equation*}
  Z(\alpha, \beta, \gamma)=
  \begin{pmatrix}
    0 & \alpha & -\overline{\alpha} & 0\\
    \overline{\alpha} & -\gamma& 0 & -\beta\\
    -\alpha & 0 & \gamma & \overline{\beta} \\
    0 & -\overline{\beta} & \beta & 0\\
  \end{pmatrix},
  \qquad \alpha, \beta \in \C, \qquad \gamma \in\R.
\end{equation*}
The objective is
\begin{equation*}
  \val(Z)
  =
  \Tr((\I_2\otimes S)Z)
  =
  \gamma.
\end{equation*}
Let us now consider the capacity constraint $Z' \coloneqq Q(c,r)-Z(\alpha, \beta, \gamma)\succeq0$.
Since $Z'_{00, 00} = Z'_{11, 11} = 0$, the matrix can only be positive semidefinite if $\alpha = \beta = 0$. The condition is thus equivalent to
\begin{equation*}
  \begin{pmatrix}
    c+\gamma&r\\
    r&c-\gamma
  \end{pmatrix}
  \succeq0.
\end{equation*}
This in turn is true if and only if
\begin{equation*}
  c+\gamma\ge0,\qquad c-\gamma\ge0,\qquad (c+\gamma)(c-\gamma)\ge r^2.
\end{equation*}
The maximum flow is therefore
\begin{equation*}
  \MaxFlow(G, K(c,r); S,T)
  =
  \sqrt{c^2-r^2}.
\end{equation*}

On the other hand, the projection-valued cut is forced:  We require $E, F \psubseteq V(G)$
form a partition of $V(G)$ such that moreover
$\mmap{E}\mmap{S} = \mmap{S}$, and $\mmap{F}\mmap{T} = \mmap{T}$. But these two conditions are equivalent to $\pproj{E} + \pproj{F} = \I_2$ and $\pproj{E}\pproj{S} = \pproj{S}$, $\pproj{F}\pproj{T} = \pproj{T}$, which is only satisfied if $\pproj{E}=\pproj{S}$ and $\pproj{F}=\pproj{T}$.
We can now write the cost of the cut in terms of $Q$ instead of $K$ as
\begin{equation*}
  \cut_K(E, F) \quad=\quad \input{diagrams/cut-cost-realigned.tikz}\enspace.
\end{equation*}
It follows that
\begin{equation*}
  \MinCut(G, K(c,r);S,T)
  =
  \Tr\bigl[(\pproj{S}\transpose \otimes \pproj{T})Q(c,r)\bigr]
  =
  c.
\end{equation*}
Consequently,
\begin{equation*}
  \MaxFlow(G,K(c,r);S,T)
  =
  \sqrt{c^2-r^2}
  \le c
  =
  \MinCut(G,K(c,r);S,T),
\end{equation*}
with equality if and only if $r=0$.

\begin{remark}
  We remark that similar gaps between other quantum notions of maximum flows
  and minimum cuts have been noted in the
  literature. An example is the max-flow min-cut problem for tensor networks
  \cite{cui2016quantum,Hastings:AsymptoticsQMaxFlowMinCut,GesLysSte:QMaxFlowBridge}. A
  tensor network is a graph whose edges represent vector spaces, and whose
  vertices may be assigned tensors matching the dimensions of the incoming
  edges. A natural capacity assignment to the edges is given by the dimensions
  of the corresponding vector spaces. This suggests to define the following
  flow problem: Any assignment of tensors to the vertices defines a linear map,
  and the flow value associated to the assignment is given by the rank of this
  map. Dually, a cut is simply a cut of the underlying graph, whose cost is
  given by the product of the edge capacities. Any cut factors the linear map
  through the vector space formed by the cut edges, so any cut upper bounds 
  the value of any flow, and hence the maximum flow is at most the minimum cut.
  On the other hand, it turns out that there exist tensor
  networks for which the minimum cut is strictly larger than the maximum flow
  \cite{cui2016quantum,GesLysSte:QMaxFlowBridge}.
  While there does not seem to be any obvious connection to our work, it would
  be interesting to investigate whether a relaxation of the cut definition could
  recover equality, just as in our case.
\end{remark}

\section{Examples}\label{sec:examples}
To better illustrate our definitions, we compute several examples of flows on
quantum graphs. On the one hand, we show that certain classical results
naturally generalise to the quantum case. For instance, the maximum flow
achievable on complete quantum graphs is given by the same formula as for the
classical complete graph. We also show that the strong phrasing of the max-flow
min-cut theorem holds for certain subclasses of quantum graphs. On the other
hand, we will see that certain non-classical phenomena arise. For instance,
capacity maps with off-diagonal elements enforce coherence between different
edges, limiting possible flow assignments. This will lead to surprising results
like the existence of connected quantum graphs that nevertheless do not admit
non-zero flows between certain source and sink sets.

\subsection{Classical Graphs Embedded in \texorpdfstring{$M_n$}{Mn}}

In \cite{spitzer2026quantum}, the authors introduce the family of quantum
graphs $X_{A, \ccdot}$. These quantum graphs are given by the embedding of the
adjacency matrix $A$ of a classical graph into $M_n$ through the diagonal
embedding. In many ways, these graphs behave exactly as their classical
counterpart, cf. \cite[Table 1]{spitzer2026quantum}, so it is natural to investigate how
flows and cuts of $X_{A, \ccdot}$ relate to flows and cuts of the classical graph $A$. A
crucial difference that arises just by viewing $A$ as a quantum graph on $M_n$
instead of a quantum graph on $\C^n$ (a classical graph) is that the projective
subsets $S$ and $T$, encoding source and sink respectively, no longer reduce to
classical subsets of the vertex set. Instead, they are in bijective
correspondence with projectors on $\C^n$. These projections may hence be
interpreted as \emph{superpositions} of sources and sinks. Another difference
is that capacities on $X_{A, \ccdot}$ have more degrees of freedom: While a
capacity on $A$ is given by an $n \times n$ matrix, a capacity on $X_{A,
\ccdot}$ is a linear map $M_n \to M_n$ and is thus represented by an $n^2
\times n^2$ matrix. We will see that the additional degrees of freedom may be interpreted as
coherence between different edges. As for sinks and sources,
classical capacities are special cases of the quantum capacities.

\begin{definition}
  Let $K$ be a capacity of a classical graph. Then $\tilde{K}\colon M_n \to
  M_n$ is the \emph{diagonal capacity}, given by
  \begin{equation*}
    \input{diagrams/diag-capac.tikz}
  \end{equation*}
\end{definition}

The name ``diagonal capacity'' is motivated by the fact that when passing to
the Choi picture by realigning, one finds that $\tilde{K}$ is precisely the
embedding of the edge capacities onto the main diagonal. We verify that
$\tilde{K}$ is indeed a valid capacity on $X_{A, \ccdot}$.

\begin{proposition}\label{prop:diagcapac}
  Let $K$ be a capacity on a classical graph $A$. Then $\tilde{K}$ is a capacity on $X_{A, \ccdot}$.
  \begin{proof}
    The adjacency operator of $X_{A, \ccdot}$ is given by
    \begin{equation*}
      \input{diagrams/xa-adjop.tikz}\enspace.
    \end{equation*}
    It thus follows that $\tilde{K}$ is supported on $X_{A, \ccdot}$:
    \begin{equation*}
      \input{diagrams/xa-capac-support.tikz}\enspace.
    \end{equation*}
    Since the diagram expressing $\tilde{K}$ in terms of $K$ is symmetric under
    180 degree rotations, we have that $\tilde{K}\transpose = \tilde{K}$ since
    $K\transpose = K$. Finally, $\tilde{K}$ is completely positive if and only
    if its realignment is a positive morphism. The realignment of $\tilde{K}$
    is given by
    \begin{equation*}
      \input{diagrams/diag-capac-realigned.tikz}\enspace,
    \end{equation*}
    but this is positive by \cref{def:cp} since $K$ is completely positive.
  \end{proof}
\end{proposition}

The diagonal capacity is precisely what is needed to simulate the flow
behaviour of a classical graph $A$ on its quantum counterpart $X_{A, \ccdot}$.

\begin{proposition}\label{prop:xaclflowstoqflows}
  Let $K$ be a capacity on a classical graph $A$ and let $S, T \psubseteq
  V(A)$. Then there is a correspondence between $S$-$T$-flows $F$ on $A$ and
  $S$-$T$-flows $\tilde{F}$ on $X_{A, \ccdot}$ such that $\val(K; F) =
  \val(\tilde{K}; \tilde{F})$.
  \begin{proof}
    Let $F$ be a flow on $A$. Just as for capacities, we define
    \begin{equation*}
      \input{diagrams/diag-flow.tikz}\enspace.
    \end{equation*}
    This is a valid flow: By the same reasoning as in \cref{prop:diagcapac},
    $\tilde{F}$ satisfies
    $\tilde{F}\transpose = -\tilde{F}$, and $\tilde{K} - \tilde{F}$ is
    completely positive. Since the embedding is symmetric under horizontal
    mirroring, we moreover have $\overline{\tilde{F}} = \tilde{F}$. Finally,
    conservation is satisfied, as the following graphical derivation
    shows.
    \begin{equation*}
      \input{diagrams/diag-flow-cons.tikz}
    \end{equation*}
    As desired, its value is given by
    \begin{equation*}
      \val(\tilde{F}) \quad=\quad \input{diagrams/diag-flow-value.tikz} \quad=\quad \val(F).
    \end{equation*}

    Conversely, suppose $\tilde{E}$ is a flow on $X_{A, \ccdot}$. Then we can define a flow on $A$ by
    \begin{equation*}
      \input{diagrams/undiag-flow.tikz}\enspace.
    \end{equation*}
    It is again obvious that $E\transpose = -E$. Now $K - E$ is completely positive since
    \begin{align*}
      \input{diagrams/undiag-flow-cp.tikz} &\quad=\quad \input{diagrams/undiag-flow-cp-2.tikz}\\
                                 &\quad=\quad \input{diagrams/undiag-flow-cp-3.tikz}\enspace.
    \end{align*}
    It satisfies conservation since
    \begin{equation*}
      \input{diagrams/undiag-flow-cons.tikz},
    \end{equation*}
    and we have
    \begin{equation*}
      \val(E) \quad=\quad \input{diagrams/undiag-flow-value.tikz} \quad=\quad\val(\tilde{E}),
    \end{equation*}
    which shows that it has the same value as the quantum flow.
  \end{proof}
\end{proposition}

Note that the correspondence is not bijective. In the Choi picture, this
corresponds to the fact that flows on $X_{A, \ccdot}$ might have off-diagonal
values, which do not matter for the value of the flow if source and sink are
classical. We now get that the following strong version of the max-flow min-cut
theorem holds for quantum graphs of the form $X_{A, \ccdot}$.

\begin{corollary}\label{cor:xastrongmaxflowmincut}
  Let $A$ be a classical graph and $S, T \psubseteq V(A)$. Then for any
  capacity $K$ on $A$, there exists a capacity $\tilde{K}$ on $X_{A, \ccdot}$
  such that we have
  \begin{align*}
    \MaxFlow(X_{A, \ccdot}, \tilde{K}; S, T) &= \MaxFlow(A, K; S, T)\\
                                   &= \MinCut(A, K; S, T)\\
                                   &= \GenMinCut(X_{A, \ccdot}, \tilde{K}; S, T) \\
                                   &= \MinCut(X_{A, \ccdot}, \tilde{K}; S, T).
  \end{align*}
  \begin{proof}
    The first equality follows from \cref{prop:xaclflowstoqflows}, the second
    equality follows from the classical max-flow min-cut theorem, and the third
    equality follows from \cref{thm:qmaxflowmincut}. To show the last
    equality, note that the value of a cut $(E, F)$ of $X_{A, \ccdot}$ is given
    by
    \begin{equation*}
      \input{diagrams/xa-cut-value.tikz}\enspace.
    \end{equation*}
    If $(E, F)$ is classical, then this is equal to
    \begin{equation*}
      \input{diagrams/xa-cut-value-classical.tikz} \quad=\quad \cut_K(E, F).
    \end{equation*}
    This immediately implies that $\MinCut(X_{A, \ccdot},
    \tilde{K}; S, T) \leq \MinCut(A, K; S, T)$. At the same time, by \cref{thm:qmaxflowmincut}
    we have $\MinCut(X_{A, \ccdot}, \tilde{K}; S, T) \geq
    \GenMinCut(X_{A, \ccdot}, \tilde{K}; S, T) = \MinCut(A, K; S, T)$.
  \end{proof}
\end{corollary}

However, there exist capacities on $X_{A, \ccdot}$ that do not correspond to
capacities on $A$. In fact, we have already seen such a capacity. The
counterexample separating the value of the maximum flow and the
projection-valued minimum cut introduced in \cref{ssec:ordcuts} is precisely the
quantum graph $X_{K_2, \ccdot}$ together with a non-classical capacity.
This also explains our observation that maximum flow / minimum generalised cut
and the ordinary minimum cut coincide for the capacity
\begin{equation*}
  K^R = \begin{bmatrix}
    c & r\\
    r & c
  \end{bmatrix}
\end{equation*}
if and only if $r = 0$: If $r = 0$, then $K$ is precisely equal to the diagonal
capacity corresponding to the classical capacity that assigns $c$ to the single
undirected edge of $K_2$. We are thus in the fully classical case, and
\cref{cor:xastrongmaxflowmincut} applies. To interpret the case where $r \neq
0$, note that the rows and columns of $K^R$ are precisely indexed by the
directed edges of $K_2$. The diagonal elements thus only impose restrictions on
edges, while off-diagonal elements impose restrictions on \emph{pairs} of
edges. Non-zero $r$ thus introduces non-trivial coherence between different
directed edges. This coherence is invisible to the projection-valued cut $(S,
T)$, but it reduces the positive slack available for an antisymmetric flow.

In summary, we find that for classical capacities and classical sources and
sinks, flows and cuts on $X_{A, \ccdot}$ behave the same as on its underlying
classical graph $A$. At the same time, $X_{A, \ccdot}$ also admits
non-classical capacities, which encode coherence between different edges; and
non-classical sources and sinks, which can be seen as superpositions of sources
and sinks. Studying general flows and cuts on $X_{A, \ccdot}$ can thus be seen
as a way to study \emph{quantum flows} and \emph{quantum cuts} of classical
graphs.

\subsection{The Complete Quantum Graph}
Usually the easiest quantum graph to consider on a quantum set $X$ is the
complete quantum graph $K_X$, which is given by the adjacency operator
\begin{equation*}
  \input{diagrams/kn-adjop.tikz}\enspace.
\end{equation*}
Curiously, in the case of flows the complete graph turns out to be the most
general case. Indeed, it is easy to see that every capacity of a (loopless)
quantum graph is also a valid capacity for $K_X$. Instead, we consider $K_X$
with unit capacities.

It turns out that flows and cuts closely mirror the classical case. Indeed,
consider how one would maximise the flow in the classical case: For any choice
$S, T \subseteq V(K_n)$ we can decompose the complete graph as 
\begin{equation*}
  \begin{tikzpicture}[
    scale=0.5,
    blob/.style = {fill=gray!25, draw=gray, line width=0.8pt},
    edge/.style = {black, line width=0.5pt}
    ]

    \foreach \ya/\yb in {1.2/2.1, 0.4/0.7, -0.4/-0.7, -1.2/-2.1}
    \draw[edge] (-4.0,\ya) -- (0,\yb+0.7);

    \foreach \ya/\yb in {2.1/1.2, 0.7/0.4, -0.7/-0.4, -2.1/-1.2}
    \draw[edge] (0,\ya+0.7) -- (4.3,\yb);

    \begin{scope}[overlay]
      \draw[edge] (-4.0,-1.0) .. controls (-2.0,-4.9) and (2.0,-4.9) .. (4.3,-1.0);
      \draw[edge] (-3.8,-0.8) .. controls (-1.6,-4.4) and (1.6,-4.4) .. (4.1,-0.8);
    \end{scope}
    \path (0,-3.95);

    \begin{scope}[shift={(-4.2,0.1)}, scale=0.85, shift={(4.2,-0.1)}]
      \draw[blob] plot[smooth cycle, tension=0.8] coordinates {
        (-5.6, 0.2) (-5.2, 1.6) (-4.1, 2.1) (-3.0, 1.5)
      (-2.8, 0.3) (-3.1,-1.2) (-4.2,-1.9) (-5.3,-1.3)};
    \end{scope}
    \node at (-4.2, 0.1) {$T$};

    \begin{scope}[yshift=0.35cm]
      \draw[blob] plot[smooth cycle, tension=0.7] coordinates {
        (-1.0, 2.9) ( 0.1, 3.2) ( 1.1, 2.8) ( 1.0, 1.2)
        ( 1.2,-0.6) ( 0.9,-2.6) (-0.2,-3.1) (-1.1,-2.7)
      (-0.9,-1.0) (-1.2, 0.9)};
      \node at (0, 0) {$K_m$};
    \end{scope}

    \begin{scope}[shift={(4.4,0.1)}, scale=0.85, shift={(-4.4,-0.1)}]
      \draw[blob] plot[smooth cycle, tension=0.8] coordinates {
        ( 3.0, 0.6) ( 3.3, 1.8) ( 4.4, 2.0) ( 5.5, 1.3)
      ( 5.7,-0.2) ( 5.1,-1.6) ( 3.9,-2.0) ( 3.1,-1.1)};
    \end{scope}
    \node at (4.4, 0.1) {$S$};
  \end{tikzpicture}
\end{equation*}
Clearly, the edges connecting $S$ to $T$ directly should have full flow. For
the remaining edges, assume that $\abs{S} \leq \abs{T}$. Then we can assign
full flow to all edges leaving $S$, at which point the central vertices all
have net in-flow $\abs{S}$. To restore flow conservation, we can then assign
flow $\abs{S}/\abs{T}$ to the edges into $T$. Clearly, the flow is maximised,
as every edge leaving $S$ has maximum flow. There are precisely $\abs{S} \cdot
(n - \abs{S})$ such edges, so taking into account the case $\abs{S} > \abs{T}$,
we have $\MaxFlow(K_n; S, T) = \min\{\abs{S} \cdot (n - \abs{S}), \abs{T} \cdot
(n - \abs{T})\}$.

Now let $X$ be a quantum set and $S, T \psubseteq X$. We can directly translate
the classical approach into the following flow map.\footnote{Note that this map should exceptionally be read from top to bottom to recover the correct flow direction. This mirrors the classical convention that an edge from $i$ to $j$ is represented by the $(i, j)$-th entry of the adjacency matrix.}
\begin{equation*}
  \input{diagrams/kn-flowmap.tikz}
\end{equation*}
Here, $m \coloneqq \min \{ \abs{S}, \abs{T} \}$. Clearly this map is real and
antisymmetric. We can also show that $K_X - F$ is completely positive. Indeed, we have
\begin{equation*}
  \input{diagrams/kn-flowmap-cp.tikz}
\end{equation*}
Except for the identity, the summands are all multiples of mutually orthogonal
projectors. In particular, it follows that the maps are all simultaneously diagonalisable,
and since $m/\abs{S}, m/\abs{T} \leq 1$, the map can be written as a diagonal
matrix with non-negative entries on the main diagonal. Clearly the resulting
morphism is positive in $\FdHilb$, and hence $K_X - F$ is completely positive.
The flow map also satisfies conservation, since
\begin{align*}
     &\input{diagrams/kn-flowmap-cons.tikz}\\
     &=\quad \input{diagrams/kn-flowmap-cons-2.tikz}\\
     &=\quad \input{diagrams/kn-flowmap-cons-3.tikz}\\
     &=\quad 0.
\end{align*}
Finally, its value is given by
\begin{align*}
  \val(F) \quad&=\quad \input{diagrams/kn-flowmap-value.tikz}\\
                       &=\quad \input{diagrams/kn-flowmap-value-2.tikz}\\
                       &=\quad \abs{S} \cdot \abs{T} + \min \{\abs{S}, \abs{T}\} \cdot \abs{R}\\
                       &=\quad \min \{\abs{S} \cdot (\abs{X} - \abs{S}), \abs{T} \cdot (\abs{X} - \abs{T}) \}.
\end{align*}

The case of cuts is analogous. Classically, the smallest cut is clearly given by the edges leaving $S$ or $T$, depending on which is smaller. We can do the same for quantum graphs by considering the cuts $(S, \compl{S})$ or $(T, \compl{T})$. Computing their values gives
\begin{equation*}
  \cut_K(S, \compl{S}) \quad=\quad \input{diagrams/kn-cut.tikz}  \enspace=\quad \abs{S} \cdot (\abs{X} - \abs{S})
\end{equation*}
and analogously $\cut_K(T, \compl{T}) = \abs{T} \cdot (\abs{X} - \abs{T})$. By always choosing the smaller cut we match the flow value we derived. By \cref{thm:qmaxflowmincut}, we have 
\begin{align*}
  \min \{\abs{S} \cdot (\abs{X} - \abs{S}), \abs{T} \cdot (\abs{X} - \abs{T})\} &\leq \MaxFlow(K_X; S, T) \\
  &\leq \MinCut(K_X; S, T)\\
  &\leq \min \{\abs{S} \cdot (\abs{X} - \abs{S}), \abs{T} \cdot (\abs{X} - \abs{T})\},
\end{align*}
so the two values must also be optimal for quantum $K_X$. We thus obtain the following result.

\begin{proposition}\label{prop:completegraph}
  Let $X$ be a quantum set. For all disjoint $S, T \psubseteq X$ it holds that $\MaxFlow(K_X; S, T) = \MinCut(K_X; S, T) = \min \{ \abs{S} \cdot (\abs{X} - \abs{S}), \abs{T} \cdot (\abs{X} - \abs{T}) \}$.
\end{proposition}

\begin{remark}
  In contrast to the symmetric quantum graph $\Gsym$, discussed in the
  following section in \cref{rem:gsymnonclassical}, the value in
  \cref{prop:completegraph} does not depend on whether the terminals are
  classical, but only on their ranks. This reflects the invariance of $K_{M_n}$
  under the full unitary group $U(n)$, which acts transitively on pairs of
  orthogonal projectors with given ranks.
\end{remark}

\subsection{Quantum Graphs with Diagonal Symmetry}

The quantum graphs we now turn to arise from a guiding principle explored in
\cite{spitzer2026quantum}: Considering some form of symmetry for quantum graphs
makes their analysis more tractable. Here, we shall consider symmetry groups
given by subgroups of the unitary group $U(n)$. Invariance under the full
unitary group leaves only the complete and the edgeless quantum graphs, whereas
the diagonal unitary group $DU(n)$ and the diagonal orthogonal group $DO(n)$
produce genuinely non-trivial families, parametrised by triples of matrices
$(A,B,C)$ with matching diagonals and denoted $X_{A,B,C}$; we refer the reader
to \cite{spitzer2026quantum} for details. Each matrix plays a distinct role:
$A$ encodes a classical graph, $C$ introduces \emph{strange edges} annotated
with a phase $\theta\in[0,2\pi)$, and $B$ is a projector attaching a purely
quantum subspace with no classical analogue. Taken together, $A$ and $C$
determine a classical model, the \emph{strange graph} $\mathfrak{G}(A,C)$,
whose edges are either classical or strange, and which can be thought of as the
classical skeleton of $X_{A,B,C}$; the genuinely quantum content lives in the
remaining matrix $B$. With this picture in mind, we now study the max flow
problem for the strange graph part of the quantum graphs $X_{A,B,C}$. We take
unit capacity, namely $K=n^2\Pi_A$, where $\Pi_A$ is the Choi support projection
of the quantum graph.

For an unordered pair $e=\{i,j\}$, $i<j$, put
$H_e=\lspan\{\ket{ij},\ket{ji}\}$.  The support block on $H_e$ has three
possibilities.  A missing edge has $\Pi_e=0$.  A classical edge has
$\Pi_e=\I_{H_e}$.  A strange edge with phase $\theta_e$ has rank-one support
$\Pi_e=\ket{u_e}\bra{u_e}$, where
\begin{equation*}
  u_e=\frac{\ket{ij}+e^{-i\theta_e}\ket{ji}}{\sqrt2}.
\end{equation*}
Equivalently, in the ordered basis $(\ket{ij},\ket{ji})$, this strange block
is

\begin{equation*}
  \frac{1}{2} \begin{bmatrix}
    1 & e^{i\theta_e}\\
    e^{-i\theta_e} & 1
  \end{bmatrix}.
\end{equation*}

Let $R_i=\ket{i}\bra{i}$.  For disjoint coordinate source and sink sets
$S_0,T_0\subseteq[n]$, set $P_S=\sum_{s\in S_0}R_s$ and
$P_T=\sum_{t\in T_0}R_t$.  Let $G_{\mathrm{cl}}$ be the ordinary graph
whose edges are precisely the rank-two classical edges of the strange graph,
each with unit capacity.  The strange edges are omitted from
$G_{\mathrm{cl}}$.

\begin{proposition}\label{prop:xacunit}
  For the unit-capacity $X_{A,\ccdot,C}$ graph above,
\begin{equation*}
  \MaxFlow(X_{A, \ccdot, C};P_S,P_T)
  =
  \MaxFlow(G_{\mathrm{cl}};S_0,T_0)
  =
  \GenMinCut(X_{A, \ccdot, C};P_S,P_T).
\end{equation*}
In particular, if every edge of the underlying ordinary strange graph is
strange, then the coordinate $S_0$-to-$T_0$ quantum maximum flow is zero,
for every choice of phases.
\end{proposition}

\begin{proof}
Let $Z$ be feasible for the quantum max-flow SDP, that is, self-adjoint,
antisymmetric, with $Q-Z\succeq0$ and conserved on $R$. For simplicity, we suppress the 
$n^2$ factor of the unit capacity $Q = n^2 X_{A, \ccdot, C}^R$.

First compress $Z$ to
a classical block $H_{\{i,j\}}$, $i<j$.  On the ordered basis
$(\ket{ij},\ket{ji})$ the transpose $(-)\transpose$ swaps the diagonal entries and
preserves the off-diagonal ones, so for a self-adjoint antisymmetric matrix the
off-diagonal entry is fixed by the transpose and must equal its own negative; it
therefore vanishes, and the compression has the form
\begin{equation*}
  Z_{\{i,j\}}
  =
  g_{ij}(\ket{ji}\bra{ji}-\ket{ij}\bra{ij})
  \quad\text{for some }g_{ij}\in\R .
\end{equation*}
On this block the
compression of $Q-Z$ is
$\operatorname{diag}(1+g_{ij},1-g_{ij})$, so positivity gives
$\abs{g_{ij}}\leq1$.  Thus every rank-two block supplies exactly one
 unit of capacity.

Next compress $Z$ to a strange edge block.  As above, self-adjointness and
antisymmetry force the compression to be diagonal, say
$Z_{\{i,j\}}=\operatorname{diag}(h_e,-h_e)$ for some $h_e\in\R$.  On this
block the capacity is the rank-one projection
\begin{equation*}
  Q_{\{i,j\}}
  =
  \ket{u_e}\bra{u_e}
  =
  \frac12
  \begin{bmatrix}
    1 & e^{i\theta_e}\\
    e^{-i\theta_e} & 1
  \end{bmatrix},
\end{equation*}
so the compression of $Q-Z$ has determinant
$(\tfrac12-h_e)(\tfrac12+h_e)-\tfrac14=-h_e^2\le0$.  Since this compression
is positive semidefinite, its determinant must be nonnegative, and thus
$h_e=0$.  A single strange edge therefore carries no diagonal orientable
flow.

Off-diagonal blocks between different edges may exist, but they do not
change the coordinate boundary value.  Indeed, for a Choi matrix unit,
\begin{equation*}
  \Tr\bigl(R_k\,\Tr_1(\ket{ab}\bra{cd})\bigr)
  =
  \delta_{ac}\delta_{bd}\delta_{bk}.
\end{equation*}
Hence the diagonal of $\Tr_1 Z$, and therefore the objective
$\Tr((\I_n\otimes P_S)Z)$, only sees diagonal directed Choi matrix units
inside a single edge block. Therefore every feasible
quantum flow determines a feasible classical flow on $G_{\mathrm{cl}}$ with
the same value.

Conversely, let $g$ be any feasible unit-capacity classical flow on
$G_{\mathrm{cl}}$, written antisymmetrically with $g_{ij}=-g_{ji}$.  Define
\begin{equation*}
  Z_g
  =
  \sum_{\{i,j\}\in E(G_{\mathrm{cl}}),\,i<j}
  g_{ij}(\ket{ji}\bra{ji}-\ket{ij}\bra{ij}).
\end{equation*}
Then $Z_g$ is self-adjoint, antisymmetric, and conservative exactly because $g$ is
conservative at nonterminal vertices.  The blockwise positivity calculation
above gives $K-Z_g\succeq0$, while each strange block is left as the positive
slack $\Pi_e$.  The objective of $Z_g$ is the classical value of $g$, proving
the reverse inequality.  The equality with $\GenMinCut$ is the
quantum max-flow/min-cut identity of \cref{thm:qmaxflowmincut}.
\end{proof}

\begin{example}[A two-edge strange path]
Let the vertices be $0,1,2$, with source $S=R_0$ and sink $T=R_2$.  Put
a strange edge of phase $\alpha$ on $\{0,1\}$ and a strange edge of phase
$\beta$ on $\{1,2\}$.  Take unit capacity
$Q=\ket{u_{01}}\bra{u_{01}}+\ket{u_{12}}\bra{u_{12}}$, where
\begin{equation*}
  u_{01}=\frac{\ket{01}+e^{-i\alpha}\ket{10}}{\sqrt2},
  \qquad
  u_{12}=\frac{\ket{12}+e^{-i\beta}\ket{21}}{\sqrt2}.
\end{equation*}
The underlying ordinary graph of the strange graph is the path $0-1-2$,
whose ordinary unit maximum flow from $0$ to $2$ is $1$.  However,
$G_{\mathrm{cl}}$ is empty, so \cref{prop:xacunit} gives
\begin{equation*}
  \MaxFlow(X_{A,\cdot,C}; R_0, R_2)
  =
  \GenMinCut(X_{A,\cdot,C}; R_0, R_2)
  =
  0
  \quad\text{for all phases }\alpha,\beta .
\end{equation*}

The projective cut is larger.  Since the only interior coordinate is
$\ket1$, an admissible projective cut has either $E=R_0$, $F=R_1+R_2$,
or $E=R_0+R_1$, $F=R_2$.  The first cut crosses $\{0,1\}$, the second
crosses $\{1,2\}$, and in both cases
$\Tr((R_j\otimes R_i)\ket{u_{ij}}\bra{u_{ij}})
=\abs{\langle ji,u_{ij}\rangle}^2=1/2$.  Thus
$\MinCut(X_{A,\cdot,C}, K; R_0, R_2)=1/2$.
\end{example}

It was shown in \cite[Proposition 3.18]{spitzer2026quantum} that the
\emph{antisymmetric quantum graph} $\Gasym$ is precisely the quantum graph
$X_{A, \ccdot, C}$ whose strange graph is the complete graph of strange edges
with phase $\pi$. We thus obtain the following corollary.

\begin{corollary}\label{cor:gasymzeroflow}
  For any classical $S, T \psubseteq M_n$, we have $\MaxFlow(\Gasym; S, T) = 0$.
\end{corollary}

This immediately raises the question whether its counterpart, the
\emph{symmetric quantum graph} $\Gsym$ admits larger flows. Given the results
in \cite[Table 1]{spitzer2026quantum}, it is natural to think of $\Gsym$ as
being denser than $\Gasym$. Nevertheless, we can show that $\Gsym$ also does
not support any flow between classical sources and sinks. In fact, we can show the
stronger result that this does not change for non-unit capacities.
To prove this, we first note that the capacity bound $Q - Z \succeq 0$ requires
the following necessary condition, which is generally easier to verify and
matches the classical intuition.

\begin{lemma}
  Let $F$ be a flow on a quantum graph $G$ with capacity $K$. Then $F$ is supported on $G$, that is $\Pi_A Z \Pi_A = Z$.
  \begin{proof}
    As always we denote $Q = n^2K^R$. The capacity bound condition in
    the Choi picture thus reads $Q - Z \succeq 0$. As noted above, under the
    assumption that $Q\transpose = Q$ and $Z\transpose = -Z$, we also have $Q +
    Z \succeq 0$. Now let $v \in \ker \Pi_A$. Then we have
    \begin{align*}
      \braket{v | Q + Z | v} &= \braket{v | Q | v}  + \braket{v | Z | v} = \braket{v | Z | v} \geq 0\\
      \braket{v | Q - Z | v} &= \braket{v | Q | v}  - \braket{v | Z | v} = - \braket{v | Z | v} \geq 0
    \end{align*}
    where the second equalities respectively follow from the fact that the
    capacity is supported on the graph. This implies that $\braket{v | Z | v} =
    0$ for all $v \in \ker \Pi_A$. It follows that the diagonal elements of $Z$
    outside of the support of $G$ are zero. Now let $u$ be arbitrary. Since $Q
    + Z \succeq 0$ and $Q - Z \succeq 0$, so must be their restriction to
    $\lspan\{u, v\}$. By our above observations, the diagonal element of $v$ is
    $0$, and thus by semidefiniteness so must be the off-diagonal elements in
    the same row and column. We thus get
    \begin{equation*}
      \braket{u | Q + Z | v} = \braket{u | Q | v}  + \braket{u | Z | v} = \braket{u | Z | v} = 0
    \end{equation*}
    where the second equality again follows from the support condition. By self-adjointness, the same holds for the roles of $u$ and $v$ exchanged, so we find that $\braket{u | Z | v} = 0$ if $u \in \ker \Pi_A$ or $v \in \ker \Pi_A$. It follows that $\Pi_A Z \Pi_A = Z$ as desired.
  \end{proof}
\end{lemma}

We can now show that $\Gsym$ admits no flows.

\begin{proposition}\label{prop:gsymzero}
  For any valid capacity $K$ and any classical $S, T \psubseteq M_n$, we have $\MaxFlow(\Gsym, K; S, T) = 0$.
  \begin{proof}
    We know that the projector associated to $\Gsym$ is given by
    \begin{equation*}
      \Pi_A = \frac{1}{2}(\I + \Sigma) - \frac{1}{n} \omega\omega\adjoint,
    \end{equation*}
    where $\omega$ is the maximally entangled state and $\Sigma$ is the swap
    map. The support condition for flows now says that $\Pi_A Z \Pi_A = Z$,
    which implies $Z = Z \Pi_A$, so
    for $i,j,k,\ell \in [n]$ with $i \neq j$ and $k \neq \ell$ we have
    \begin{equation*}
      Z_{ij, k\ell} = \frac{1}{2}(Z_{ij, k\ell} + Z_{ij, \ell k})
    \end{equation*}
    and hence $Z_{ij, k\ell} = Z_{ij, \ell k}$. Similarly, we have $Z = \Pi_A Z$ and hence
    $Z_{ij, k \ell} = Z_{ji, k \ell} = Z_{ji, \ell k}$. But by antisymmetry, we also have
    for diagonal elements
    \begin{equation*}
      Z_{ji, ji} = Z\transpose_{ij, ij} = - Z_{ij, ij}.
    \end{equation*}
    Together this implies that $Z_{ij, ij} = 0$ for all $i \neq j$. The elements $Z_{ii, ii}$ moreover must be $0$ by antisymmetry alone. It follows that for any valid flow $Z$ on $\Gsym$, we have
    \begin{equation*}
      \val(F) = \Tr [(\I \otimes \pproj{S})Z] = \sum_{i}\sum_{j \in S} Z_{ij, ij} = 0,
    \end{equation*}
    as claimed.
  \end{proof}
\end{proposition}

\begin{remark}\label{rem:gsymnonclassical}
  \cref{prop:gsymzero} applies only to \emph{classical} source and sink sets. It is possible to have non-zero flow between non-classical terminals. To illustrate this, we consider the case $n=2$ and let
  \begin{equation*}
    v_1=\frac{\ket{00}-\ket{11}}{\sqrt2},\qquad
    v_2=\frac{\ket{01}+\ket{10}}{\sqrt2},\qquad
    \Pi_A=\ket{v_1}\bra{v_1}+\ket{v_2}\bra{v_2},
  \end{equation*}
  so that $\Pi_A$ is the projector onto the edge space of $\Gsym$ on $M_2$. We consider the capacity $K = n^{-2} \Pi_A^\bR$, such that $Q = \Pi_A$, and consider the non-classical terminal pair
  \begin{equation*}
    p=\frac{\ket0+i\ket1}{\sqrt2},\qquad
    \pproj{S}=\ket{p}\bra{p},\qquad
    \pproj{T}=\I_2-\pproj{S}.
  \end{equation*}
  This way, the interior projector $\pproj{R}=\I_2-\pproj{S}-\pproj{T}$ vanishes and conservation imposes no constraint. Consider
  \begin{equation*}
    Z=-i\ket{v_1}\bra{v_2}+i\ket{v_2}\bra{v_1}.
  \end{equation*}
  Then $Z$ is self-adjoint; since
  $(\ket{v_1}\bra{v_2})\transpose=\ket{v_2}\bra{v_1}$, it is antisymmetric,
  $Z\transpose=-Z$; it is supported on $\Gsym$ because its range lies in
  $\lspan\{v_1,v_2\}=\img \Pi_A$; and $\Pi_A-Z\succeq0$, since in the
  $(v_1,v_2)$ basis
  $Z=\bigl(\begin{smallmatrix}0&-i\\i&0\end{smallmatrix}\bigr)$ has eigenvalues
  $\pm 1$, so $\Pi_A-Z$ has eigenvalues $0$ and $2$ on $\img \Pi_A$ and vanishes
  on its complement. Thus $Z$ is feasible for the flow SDP.

  It remains to compute the value of this flow. In the Choi basis $(\ket{00},\ket{01},\ket{10},\ket{11})$ we have
  \begin{equation*}
    Z=\frac12
    \begin{bmatrix}
      0 & -i & -i & 0\\
      i & 0 & 0 & -i\\
      i & 0 & 0 & -i\\
      0 & i & i & 0
    \end{bmatrix},
    \qquad
    \I_2\otimes\pproj{S}=\frac12
    \begin{bmatrix}
      1 & -i & 0 & 0\\
      i & 1 & 0 & 0\\
      0 & 0 & 1 & -i\\
      0 & 0 & i & 1
    \end{bmatrix},
  \end{equation*}
  and a direct multiplication shows that $(\I_2\otimes\pproj{S})Z$ has diagonal entries $(\tfrac14,\tfrac14,\tfrac14,\tfrac14)$. Hence
  \begin{equation*}
    \val(F)=\frac{1}{4}\Tr\bigl[(\I_2\otimes\pproj{S})Z\bigr]= \frac{1}{4},
  \end{equation*}
  so $\MaxFlow(\Gsym, \pproj{S},\pproj{T}) > 0$. The mechanism is that the only non-zero matrix elements of $Z$ are the purely imaginary coherences $Z_{i\ell,ij}$ with $\ell\neq j$, between different Choi units; a diagonal projector sees only the diagonal entries $Z_{ij,ij}=0$, whereas the imaginary phase of $p$ picks up the coherences. Indeed, if $p$ were replaced by the real vector $(\ket0+\ket1)/\sqrt2$ (also non-diagonal), the same $Z$ would have value $0$. 
\end{remark}

\section{Conclusion}\label{sec:conclusion}
In this work, we have introduced a framework for flows and cuts of quantum
graphs. Our central definitions of capacity and flow maps naturally arise from
the classical concepts of capacity and flow matrices through the standard
non-commutative dictionary. Dually, we defined cuts in terms of projective
subsets, which have proven in \cite{spitzer_unification_2026} to be a natural
building block in defining properties of quantum graphs. We have shown that
under these definitions, quantum graphs satisfy a \emph{max-flow min-cut
theorem}. The theorem turns out to be slightly weaker in the quantum case
compared to the classical case. Indeed, we found that the maximum flow
achievable on a quantum graph is equal to the minimum \emph{generalised cut},
which is the quantum analogue of the continuous relaxation of the ordinary cut.

Furthermore, we have computed cuts and flows for a number of example quantum
graphs, starting with \emph{embeddings} of classical graphs into other quantum
sets. We have shown that up to a rescaling depending on the ambient quantum
set, flows and cuts are invariant under this embedding and reduce to the
classical case. In particular, we recover the strong version of the max-flow
min-cut theorem for this family of quantum graphs. Another such family turned
out to be the complete graphs on arbitrary vertex sets. We were able to show
that their maximum flows and minimum cuts are given by the same formula in
terms of the cardinality of the terminals and their complements, regardless of
the quantum set. The benefit of considering different quantum vertex sets lies in the
choice of valid source and sink sets. For example, a classical graph only
allows classical subsets as terminals, while embedding it into $M_n$ allows
arbitrary superpositions of classical terminals. This provides a route towards
studying quantum flows of classical graphs.

We concluded by considering a subclass of the ABC graphs of
\cite{spitzer2026quantum}. In particular, we have seen that flows on $X_{A,
\ccdot, C}$ between classical terminals are determined by the underlying
strange graph. Concretely, classical edges admit the same flow as their
classical counterparts, while strange edges admit no flow. This seems to
suggest that flow between classical terminals requires classical graph
structure. A natural next question is how the purely quantum contribution of
the $B$ projector of $X_{A, B, C}$ influences this setting. Since $B$ is not
encoded in the strange graph, it could for instance be used to relax the flow
conservation at vertices.

Another interesting question is whether \emph{Menger's theorem} holds for
quantum graphs. The theorem says that a graph is $k$-connected if and only if
there exist at least $k$ vertex-disjoint paths between any two vertices.
Classically, this is a corollary of the max-flow min-cut theorem. Quantumly,
there are two obstructions. The first is that quantumly, the equality between
the maximum flow and the minimum cut becomes an inequality. However, this
applies to the general setting of arbitrary capacities, while the setting of
Menger's theorem is more restricted. The second restriction is that so far
there exists no quantum analogue of paths, and it is questionable whether such
a concept is meaningful. Some of our results seem to suggest that Menger's
theorem would fail for any reasonable definition of paths. Concretely, we have
shown in \cref{cor:gasymzeroflow} that there exist connected quantum graphs
that admit no non-zero flows between certain source and sink sets.
Nevertheless, it would be interesting to investigate why it fails and in which
cases it holds.

\bigskip

\noindent\textbf{Acknowledgements.} This research was supported by the ANR project \href{https://www.ceremade.dauphine.fr/dokuwiki/anr-tagada:start}{TAGADA}. The first author moreover acknowledges support from the CNRS 80Prime grant ``QuantGraphe''.

\medskip

\noindent\textbf{Author contributions.} The authors contributed equally to the research. GLS did the majority of the writing. Both authors reviewed and approved the final version.

\medskip

\noindent\textbf{AI Statement.} AI was used for proofreading. No AI was used for research or for the preparation of the initial manuscript. 

\printbibliography

\appendix

\section{Cut-Flow Duality for Arbitrary Quantum Sets}\label{app:cutflowgeneral}

We claimed that the max-flow min-cut theorem proved in \cref{sec:mn} for quantum graphs on $M_n$ can be extended to graphs on arbitrary quantum sets. We prove this here. Let us briefly recall what we have proved so far: We have shown that by considering the Choi matrices of the capacity and flow maps $K, F\colon M_n \to M_n$, the definitions of cuts and flows become statements about positive-semidefinite matrices on $M_n$. This allowed us to phrase both the maximum flow and the minimum generalised cut as the optima of two mutually dual SDPs. This procedure becomes much more involved when considering quantum graphs $G$ with $V(G) \neq M_n$. In principle, we could proceed along the same lines. However, as a finite-dimensional \Cstar-algebra, $V(G)$ is a direct sum of matrix algebras, so computing the Choi matrices will introduce cross-terms that make both stating the SDPs and proving their duality messy. Moreover, we translated the definitions of cuts and flows into the Choi picture by means of graphical derivations, which we cannot perform in the general case without additional machinery.

\subsection*{SDPs on Arbitrary Quantum Sets}

The main technical tool that addresses both of these obstacles is the bicategorical formulation of $\dagger$-SSFMs in $\FdHilb$ in terms of $\dagger$-bimodules in $\TFdHilb$. We will not repeat the definitions here and instead refer the reader to \cite{spitzer_unification_2026} or \cite{verdon_covariant_2024}.
The crucial observation is that we can perform the diagrammatic derivations of the properties in the Choi picture in the corresponding modified graphical calculus. Firstly, for an arbitrary quantum set $X = \bigoplus_{i = 1}^k M_{n_i}$, a map $F \colon X \to X$ can be expressed in $\TFdHilb$ as
\begin{equation*}
  \input{diagrams/x-f.tikz}
\end{equation*}
where the shaded area corresponds to $X' \cong \C^k$, and the wires are given by $Y \coloneqq \bigoplus_{i = 1}^k \C^{n_i}$. The realignment procedure\footnote{Note that we include a normalisation by an invertible $2$-morphism, analogous to how we normalised the realignment on $M_n$. This differs from the literature, in particular from \cite{verdon_covariant_2022}, but it will be easy to see that all the existing results we use generalise.} then yields the map
\begin{equation*}
  \input{diagrams/x-f-realigned.tikz}
\end{equation*}
which is a linear map $\overline{Y} \otimes Y \to \overline{Y} \otimes Y$ which is additionally a $(\C^k, \C^k)$-bimodule homomorphism. Clearly, this realignment is an isomorphism, and we can recover the original $X \to X$ map by bending the wires back and undoing the scaling. It can then be shown that Choi's theorem holds in this setting, that is $F$ is completely positive if and only if $F^R$ is positive semidefinite \cite[Theorem 4.13]{verdon_covariant_2022}.

We can now derive the definitions of flows and cuts in the Choi picture using the same diagrams as in the $M_n$ case, only now including shaded regions. We start with the value of a flow. Let $F\colon X \to X$ be a flow. Then we have
\begin{equation*}
  \val(F) \quad=\quad \input{diagrams/flow-value-2realigned.tikz}
\end{equation*}
Just as for the $M_n$ case, we absorb the normalisation into a new variable, letting
\begin{equation*}
  \input{diagrams/sdp-z-x.tikz} \quad\eqqcolon\quad D_{n^2}(F^R).
\end{equation*}
The value of the flow thus becomes
\begin{equation*}
  \val(F) \quad=\quad \input{diagrams/flow-value-2realigned-z.tikz} \quad=\quad \Tr[(I \otimes \pproj{S})Z],
\end{equation*}
where the second equality follows from the fact that the duality morphisms of $Y$ are given by
\begin{equation*}
  \eta\colon \C \to \overline{Y} \otimes_{\C^k} Y, 1 \mapsto \bigoplus_{i = 1}^k \sum_{j = 1}^{n_i} \ket{j} \otimes_{\C^k} \ket{j}
\end{equation*}
and its adjoint, so the shaded trace reduces to the standard trace. Reality and antisymmetry of $F$ translate straightforwardly  to self-adjointness and antisymmetry of $Z$, fully analogous to the $M_n$ case. In particular, the argument in \cref{lem:trposerealigncommute} goes through without modification in the presence of shaded regions. The flow conservation condition becomes
\begin{equation*}
  \input{diagrams/conservation-2realigned.tikz} \quad\stackrel{!}{=}\quad 0,
\end{equation*}
which is equivalent to $\pproj{R} \Tr_1 Z = 0$. Finally, by the generalised Choi theorem mentioned above, the capacity constraint that $K - F$ be completely positive becomes $Q - Z \succeq 0$, where $Q = D_{n^2}(K^R)$. The duality between the Schur product and composition witnessed by the realignment continues to hold in the $2$-categorical setting, so the support condition on $K$ is still equivalent to $\Pi_A Q \Pi_A = Q$, and by the same arguments as above, symmetry of $K$ is equivalent to symmetry of $Q$, and complete positivity of $K$ is equivalent to positive-semidefiniteness of $Q$.

Now $\overline{Y} \otimes Y \cong M_N$, where $N \coloneqq \sum_{i = 1}^k n_i$, so we can see $Z$ simply as a map $M_N \to M_N$ by forgetting the bimodule structure. This means that we can write the max-flow problem on the vertex set $X$ as an ordinary SDP over $M_N$. The bimodule constraint can be encoded as commutation with the left and right multiplication maps $L_x$ and $R_x$ for all $x \in X'$. Concretely, we obtain the following optimisation problem:

\begin{equation}\label{eq:PP}
\begin{array}{ll}
\text{maximize} &
\lrangle{\I \otimes \pproj{S}, Z} \\
\text{subject to} &
  Z=Z\adjoint, \\
& Z\transpose=-Z, \\
& \pproj{R}\operatorname{Tr}_1 Z=0, \\
& Q-Z\succeq 0 \\
& [Z, L_x] = 0 ~\forall x \in X' \\
& [Z, R_x] = 0 ~\forall x \in X'
\end{array}
\tag{P''}
\end{equation}

By linearity, we can replace the universal quantification by the corresponding statements about basis vectors. This yields a genuine SDP. Moreover, up to the bimodule constraints, it is precisely equal to (\ref{eq:P})! The difference is in the interpretation of the variables: Suppose we are given a capacity $K$. In the $M_n$ case, we set $Q = n^2 K^R$, while in the $X$ case, we set $Q = D_{n^2}(K^R)$. Similarly, if we want to recover the flow map $F$ from the optimiser, we interpret $Z = n^2 F^R$ and $Z = D_{n^2}(F^R)$ respectively.

\smallskip

The case of the dual program is analogous. It turns out that the graphical argument in \cref{eq:dualcutinterp} goes through in the shaded calculus, replacing the scalar normalisation of $K^R$ by $Q = D_{n^2}(K^R)$. It follows that the minimum generalised cut of a graph with quantum vertex set $X$ is captured by the following SDP, which again only differs from (\ref{eq:DD}) by the additional bimodule constraints.

\begin{equation}\label{eq:DDD}
  \begin{array}{ll}
    \text{minimize} &
    \lrangle{\I \otimes \pproj{W} - \Lambda_1, Q} \\
    \text{subject to} & \Lambda_1\adjoint = \Lambda_1,\\
                      & \Lambda_1\transpose = \Lambda_1,\\
                      & \pproj{W}\adjoint = \pproj{W}, \\
                      & (\pproj{S} + \pproj{T})\pproj{W}(\pproj{S} + \pproj{T}) = \pproj{S},\\
                      & \I \otimes \pproj{W} - \Lambda_1 \succeq 0\\
                      & [\pproj{W}, R_x] = 0 ~\forall x \in X'\\
                      & [\Lambda_1, L_x] = 0 ~\forall x \in X' \\
                      & [\Lambda_1, R_x] = 0 ~\forall x \in X'
  \end{array}
  \tag{D''}
\end{equation}

\subsection*{The Proof}

We are now ready to prove the max-flow min-cut theorem for arbitrary quantum sets. We restate it here for the reader's convenience.

\thmmfmc*

To prove the theorem in full generality, we show that the bimodule constraints of (\ref{eq:PP}) and (\ref{eq:DDD}) can be eliminated. Concretely, if $S, T, R \psubseteq X$ and $Q$ is a bimodule homomorphism, then the optima of (\ref{eq:P}) and (\ref{eq:DD}) are attained at bimodule homomorphisms. Consequently, the values of (\ref{eq:P}) and (\ref{eq:PP}) as well as (\ref{eq:DD}) and (\ref{eq:DDD}) coincide. Then by the result for $M_n$, the values of (\ref{eq:PP}) and (\ref{eq:DDD}) must also coincide. In particular, the value of a maximum flow is equal to that of a minimum generalised cut.

\begin{remark}
In fact, this approach shows a slightly stronger result: As was noted in \cite[p. 27]{spitzer_unification_2026}, any quantum graph on a quantum set $X$ may be interpreted as a quantum graph on $Y$ if $X' \supseteq Y'$. For instance, this is the relationship between a classical graph $A$ on $\C^n$ and the graph $X_{A, \ccdot}$ on $M_n$. By showing that the SDPs can be freely restricted to arbitrary quantum sets, we show that up to global scaling the max-flow and min-cut values are invariant under this reinterpretation. In particular, this implies that flows on $A$ coincide with classical flows on $X_{A, \ccdot}$, which we have shown explicitly in \cref{prop:xaclflowstoqflows}.
\end{remark}

For the proof we consider the twirling map 
\begin{equation*}
  \mathcal{T}\colon B(\overline{Y} \otimes Y) \to B(\overline{Y} \otimes Y),\quad F \mapsto \int_{U(X')} \int_{U(X')} (\overline{u} \otimes v) F (\overline{u} \otimes v)\adjoint \diff u \diff v,
\end{equation*}
where $U(X')$ is the unitary group of $X'$, represented as blockwise multiples of the identity, and integration is over the Haar measure. This is the orthogonal projector onto the space of $(X', X')$-bimodule homomorphisms: The tensor product $\overline{Y} \otimes Y$ splits into blocks $\overline{\C^{n_i}} \otimes \C^{n_j}$ with corresponding projections $P_{ij}$. Writing $u = \bigoplus_i \lambda_i I_{n_i}$ and $v = \bigoplus_j \mu_j I_{n_j}$ with $|\lambda_i| = |\mu_j| = 1$, the unitary $\overline{u} \otimes v$ acts on the $(i,j)$-th block as the phase $\overline{\lambda_i}\mu_j$. Hence conjugation by $\overline{u} \otimes v$ multiplies the block $P_{ij} F P_{i'j'}$ of $F$ by $\overline{\lambda_i}\mu_j\lambda_{i'}\overline{\mu_{j'}}$. Since the phases are independent and uniformly distributed on the unit circle,
\begin{equation*}
  \mathbb{E}\big[\overline{\lambda_i}\lambda_{i'}\,\mu_j\overline{\mu_{j'}}\big] = \mathbb{E}\big[\overline{\lambda_i}\lambda_{i'}\big]\,\mathbb{E}\big[\mu_j\overline{\mu_{j'}}\big] = \delta_{ii'}\delta_{jj'},
\end{equation*}
so averaging removes all off-diagonal blocks, and we get
\begin{equation*}
  \mathcal{T}(F) = \sum_{i,j=1}^k P_{ij} F P_{ij}.
\end{equation*}
Hence $\mathcal{T}$ is the orthogonal projection onto the block-diagonal matrices. These are exactly the matrices that correspond to $(X', X')$-bimodule homomorphisms and hence satisfy the commutation constraints of (\ref{eq:PP}). Moreover, $\mathcal{T}$ leaves these maps unchanged, so it is surjective onto this space.

\paragraph{The Flow SDP.} Crucially, $\mathcal{T}$ preserves the objective of (\ref{eq:P}). Indeed, if $S \psubseteq X$ then $\pproj{S}$ is a right $X'$-module homomorphism, and hence $\identity \otimes \pproj{S}$ is a $(X', X')$-bimodule homomorphism. We thus get
\begin{equation*}
  \langle \I \otimes \pproj{S}, \mathcal{T}(Z) \rangle = \langle \mathcal{T}(\I \otimes \pproj{S}), Z \rangle = \langle \I \otimes \pproj{S}, Z \rangle
\end{equation*}
as desired. Moreover, $\mathcal{T}$ preserves the constraints. This is clear for the constraint $Z = Z\adjoint$. We further have
\begin{align*}
  \mathcal{T}(Z)\transpose \enspace&=\enspace \int_{U(X')}\int_{U(X')}\enspace \input{diagrams/twirling-pres-transpose.tikz} \enspace \diff u \diff v 
                                   \enspace=\enspace \int_{U(X')}\int_{U(X')} \input{diagrams/twirling-pres-transpose-2.tikz} \diff u \diff v \\
                                   &=\enspace - \int_{U(X')}\int_{U(X')} \enspace \input{diagrams/twirling-pres-transpose-3.tikz}\enspace \diff u \diff v \enspace=\enspace -\mathcal{T}(Z).
\end{align*}
The conservation condition works out as
\begin{align*}
  \pproj{R} \Tr_1 \mathcal{T}(Z) \enspace&=\enspace \int_{U(X')}\int_{U(X')} \enspace \input{diagrams/twirling-pres-cons.tikz} \enspace \diff u \diff v \enspace=\enspace \int_{U(X')} \int_{U(X')} \enspace \input{diagrams/twirling-pres-cons-2.tikz} \enspace\diff u \diff v \\
  \enspace&=\enspace \int_{U(X')} \enspace \input{diagrams/twirling-pres-cons-3.tikz} \enspace \diff v \enspace=\enspace 0,
\end{align*} 
where the second equality follows from the fact that $\pproj{R}$ is a right $X'$-module homomorphism and hence commutes with all $v \in U(X')$.
Finally, if $Q$ is the realignment of a capacity on a quantum graph with vertex set $X$, then it must also be a $(X', X')$-bimodule homomorphism. We thus have
\begin{equation*}
  Q - \mathcal{T}(Z) \succeq 0 \quad\IFF\quad \mathcal{T}(Q - Z) \succeq 0
\end{equation*}
and it is obvious that $\mathcal{T}$ preserves positivity.
It follows that $\mathcal{T}$ is a surjective map from the feasible region of (\ref{eq:P}) into the feasible region of (\ref{eq:PP}) which preserves the objective. Consequently, whenever $Z$ is an optimiser for the max-flow SDP then so is $\mathcal{T}(Z)$, and $D_{n^{-2}}(\mathcal{T}(Z)^\bR) \colon X \to X$ is a maximum flow. Moreover, (\ref{eq:P}) and (\ref{eq:PP}) have the same value.

\paragraph{The Cut SDP.} We proceed analogously for (\ref{eq:DDD}). A slight difference arises through the fact that we are now optimising over two variables: the weight function $\pproj{W}$ itself and the symmetric correction term $\Lambda_1$. A valid optimiser for the quantum set $X$ requires $\Lambda_1$ to be a $(X', X')$-bimodule homomorphism, while $\pproj{W}$ must be a right $X'$-module homomorphism. For that reason, we introduce a second twirling map $\mathcal{T}'$, defined as
\begin{equation*}
  \mathcal{T'} \colon B(Y) \to B(Y), F \mapsto \int_{U(X')} u F u\adjoint \diff u.
\end{equation*}
This is again a standard definition and yields the orthogonal projector onto the space of $(X', X')$-bimodule homomorphisms on $Y$. Now the action of $X'$ on $Y$ is given by blockwise scalar multiplication, and since scalar multiplication is commutative, bimodules and right-modules coincide. $\mathcal{T}'$ is thus also the orthogonal projector onto the space of right $X'$-module homomorphisms on $Y$. 

Again, we claim that if $\pproj{W}$ and $\Lambda_1$ are optimisers of the unconstrained SDP (\ref{eq:DD}), then so are $\mathcal{T'}(\pproj{W})$ and $\mathcal{T}(\Lambda_1)$. Indeed, the objective is preserved: We have
\begin{equation}\label{eq:tprimeext}
  \identity \otimes \mathcal{T'}(\pproj{W}) = \identity \otimes \int_{U(X')} v \pproj{W} v\adjoint \diff v = \int_{U(X')} \int_{U(X')} \overline{u} \identity \overline{u}\adjoint \otimes v \pproj{W} v\adjoint \diff u \diff v = \mathcal{T}(\identity \otimes \pproj{W}),
\end{equation}
and hence
\begin{equation*}
  \langle \identity \otimes \mathcal{T'}(\pproj{W}) - \mathcal{T}(\Lambda_1), Q \rangle = \langle \mathcal{T}(I \otimes \pproj{W} - \Lambda_1), Q \rangle = \langle \identity \otimes \pproj{W} - \Lambda_1, \mathcal{T}(Q) \rangle = \langle \identity \otimes \pproj{W} - \Lambda_1, Q \rangle
\end{equation*}
as desired. The constraints are also preserved: The self-adjointness and symmetry constraints are preserved by the same reasoning as above. The positivity condition again follows from positivity of $\mathcal{T}$, using \cref{eq:tprimeext}. It remains to show that the condition $(\pproj{S} + \pproj{T})\pproj{W}(\pproj{S} + \pproj{T}) = \pproj{S}$ is preserved. Since $S, T \psubseteq X$, it holds that $(\pproj{S} + \pproj{T})$ is a right $X'$-module homomorphism and thus commutes with all $u \in U(X')$. We thus get
\begin{equation*}
  (\pproj{S} + \pproj{T})\mathcal{T}'(\pproj{W})(\pproj{S} + \pproj{T}) = \int_{U(X')} u (\pproj{S} + \pproj{T})\pproj{W}(\pproj{S} + \pproj{T}) u\adjoint \diff u = \mathcal{T}'((\pproj{S} + \pproj{T})\pproj{W}(\pproj{S} + \pproj{T})) = \mathcal{T'}(\pproj{S}) = \pproj{S},
\end{equation*}
as desired. This shows that also (\ref{eq:DD}) and (\ref{eq:DDD}) have the same value. Consequently (\ref{eq:PP}) and (\ref{eq:DDD}) have the same value, and we get $\MaxFlow(G, K; S, T) = \GenMinCut(G, K; S, T)$ for all quantum sets $X$, all disjoint $S, T \psubseteq X$, all quantum graphs $G$ on $X$, and all capacities $K$ on $G$. 

\smallskip

Finally, the second part of the theorem---the upper bound by the minimum ordinary cut---follows from exactly the same argument as in \cref{ssec:ordcuts} since the constructed feasible point is feasible for (\ref{eq:DDD}) as well.

\end{document}